\documentclass[11pt,reqno]{amsart}

\usepackage[T1]{fontenc}
\usepackage{graphicx} 
\usepackage[margin= 1.0 in]{geometry}
\usepackage{amsthm}
\usepackage{amsmath}
\usepackage{mathrsfs}

\usepackage{amssymb}
\usepackage{libertinus}

\usepackage{xtab}
\usepackage{tikz-cd} 
\usepackage{ytableau}
\usepackage{xfrac}
\usepackage{stackrel}
\usepackage{color}
\usepackage{enumerate}
\usepackage{mathtools}
\usepackage{pgfplots}
\pgfplotsset{compat=1.18}
\usepackage{float} 
\usepackage{enumitem}

\usepackage{stmaryrd}

\usepackage{xcolor}
\usepackage{todonotes}

\usepackage{aliascnt}
\usepackage[bookmarks=true, bookmarksdepth=2]{hyperref}
\usepackage[nameinlink,capitalise]{cleveref}
\hypersetup{
   colorlinks=true,
   citecolor=blue,
   linkcolor=violet
   }

\newtheorem{theorem}{Theorem}[section]
\newtheorem{corollary}[theorem]{Corollary}
\newtheorem{proposition}[theorem]{Proposition}

\newaliascnt{lemma}{theorem}
\newtheorem{lemma}[lemma]{Lemma}
\aliascntresetthe{lemma}

\theoremstyle{definition}
\newtheorem{remark}[theorem]{Remark}

\newtheorem{definition}[theorem]{Definition}
\newtheorem{conjecture}[theorem]{Conjecture}
\newtheorem{notation}[theorem]{Notation}

\crefname{inequality}{inequality}{inequalities} 
\Crefname{inequality}{Inequality}{Inequalities} 
\crefname{lemma}{Lemma}{Lemmas}
\Crefname{lemma}{Lemma}{Lemmas}

\numberwithin{equation}{section}

\def\Coh{\ensuremath{\mathrm{Coh}}}
\def\Hom{\ensuremath{\mathrm{Hom}}}
\def\RHom{\mathbf{R}\mathscr{H}om}
\def\Db{\ensuremath{\mathrm{D}^b}}
\def\R{\ensuremath{\mathbb R}}
\def\K{\ensuremath{\mathrm{K}}}

\def\Stab{\ensuremath{\mathrm{Stab}}}
\newcommand{\im}{\operatorname{Im}}
\newcommand{\re}{\operatorname{Re}}
\newcommand{\ch}[1]{\operatorname{ch}_{{#1}}}

\newcommand{\rk}{\operatorname{rk}}
\newcommand{\td}{\operatorname{td}}
\newcommand{\dist}{\operatorname{dist}}

\makeatletter
\renewcommand\subsection{%
  \@startsection{subsection}{2}{\z@}%
    {3.0ex \@plus 1ex \@minus .2ex}
    {1.5ex \@plus .2ex}
    {\normalfont\bfseries}
}
\makeatother

\title{Bridgeland stability and Gieseker stability I: Limit stability}
\author{Chunkai Xu}
\dedicatory{To mathematics as we have known it.}

\begin{document}
\begin{abstract}
   We study the large volume limit of Bridgeland stability conditions on smooth projective varieties. For each fixed ray in a suitable angular sector, we prove that an object of $\Db(X)$ is Bridgeland (semi)stable in the corresponding hearts for all sufficiently large radii if and only if it is a pure $H$-Gieseker (semi)stable sheaf. In this sector, we identify the category of coherent sheaves as the limiting heart in the sense of eventual membership. We also identify a perverse coherent limiting heart in the opposite angular sector.
   
   The key ingredients are uniform HN-phase estimates for bounded families of pure sheaves, coefficient bounds for the coherent images of potential destabilizers, and a criterion forcing narrow exact sequences in the varying hearts to be bounded saturated coherent sequences.
\end{abstract}

\address{Mathematics Institute \\
 University of Warwick \\
 CV4 7AL \\
 United Kingdom}

\email{Chunkai.Xu@warwick.ac.uk}

\keywords{Bridgeland stability conditions; Gieseker stability; large volume limit; limits of hearts; $t$-structures; pure coherent sheaves; perverse coherent sheaves.}
\subjclass[2020]{Primary 14F08; Secondary 14F06, 18G80}

\maketitle
\tableofcontents

\section{Introduction}

Bridgeland stability conditions were introduced in \cite{bridgeland_stability_2007} as a notion of stability for objects in a triangulated category. Recently, based on Yucheng Liu's construction \cite{liu21productstability} of stability conditions on products of elliptic curves, Li \cite{li2026remarkstabilityconditionssmooth} constructs stability conditions on arbitrary smooth projective varieties by passing through projective space.

For stability conditions on smooth projective varieties, one of the basic questions is how they recover Gieseker stability at the large volume limit. On surfaces, the large volume comparison between Bridgeland stability and Gieseker stability is already visible in Bridgeland's work on K3 surfaces \cite{bridgeland_stability_k3}. Bayer’s polynomial stability conditions include Simpson–Gieseker stability and provide a framework for large-volume limits; on surfaces, he proves a fixed-object comparison with genuine Bridgeland stability conditions; see \cite[§2.1 and Proposition 4.1]{Bayer09PolynomialStabilityCondition}.

In this paper, we prove a comparison along fixed rays, with uniform thresholds for each numerical class, on any smooth projective variety.

\subsection{Main theorem}
Fix a projective embedding $i:X\hookrightarrow\mathbb{P}^N$ and the induced polarization $H=i^*\mathcal O_{\mathbb P^N}(1)$. The stability conditions considered here have central charge
\[
 Z_{z,X}(E)=-\int_X e^{-zH}\ch{}(E)\frac{\td(X)}{i^*\td(\mathbb P^N)},
 \qquad z=re^{i\theta},
\]
and skyscraper sheaves have phase $1$. For each fixed $\theta\in(0,\pi)$ they are defined for all sufficiently large $r$; see \Cref{StabilityConditionOnSmoothProjectiveVarieties}. Write
\[
 \sigma_{r,\theta,X}=(Z_{r,\theta,X},\mathscr P_{r,\theta,X}),
 \qquad \mathscr A_{r,\theta,X}=\mathscr P_{r,\theta,X}((0,1]).
\]

\begin{theorem}\label{Thm:MainComparison}
    Fix $\theta\in (\pi-\pi/N,\pi)$ and $v\in K_{\mathrm{num}}(X)$. There is $R = R(v,\theta) > 0$, such that for any object $E\in\Db(X)$ of class $v$,
    \[
    E \text{ is } \sigma_{r,\theta, X} \text{-(semi)stable in } \mathscr{A}_{r,\theta, X} \text{ for all } r\geq R
    \quad\Longleftrightarrow\quad
    E \text{ is a pure $H$-Gieseker (semi)stable sheaf}.
    \]
\end{theorem}
Throughout, the embedding and polarization are fixed, and their contribution to the threshold is suppressed in the notation $R(v,\theta)$. The threshold is common to all objects of the fixed numerical class, but the condition on the left is stability for every radius in a ray tail. 

The implication from Gieseker stability uses boundedness of the Gieseker semistable sheaves of that class; the implication to Gieseker stability tests each object along its entire tail. A comparison at each individual large radius is a stronger assertion and is not proved here; it will be addressed in the sequel \emph{Bridgeland stability and Gieseker stability II: The Gieseker chamber}.

The implication from eventual Bridgeland stability to Gieseker stability begins by identifying which fixed objects eventually lie in the moving hearts. The embedding and angle in the following proposition are those of \Cref{Thm:MainComparison}.

\begin{proposition}[\Cref{CoherentSheavesInLargeVolumeHeart}, \Cref{LimitImpliesCoherent}, \Cref{Cor:UniformHeartOnX}]\label{Intro:LimitOfHeart}
\[
    A\in \Db(X) \text{ lies in } \mathscr{A}_{r,\theta,X} \text{ for all } r\gg0
    \quad\Longleftrightarrow\quad
    A\in\Coh(X).
\]
\end{proposition}
The 'coherence implies heart membership' assertion is uniform on bounded families of coherent sheaves, while the converse concerns each fixed derived object. Thus every fixed coherent short exact sequence becomes exact in the moving heart for all sufficiently large radii. To test eventual stability, it is enough to apply this to one chosen subsheaf and its quotient.

The leading nonzero term of
\[
 \im\bigl(Z_{r,\theta,X}(F)\overline{Z_{r,\theta,X}(E)}\bigr)
\]
then compares the reduced Hilbert polynomials of sheaves of the same dimension. This proves the implication from limit stability to Gieseker stability, including purity and the stable case; see \Cref{LimitsImpliesGieseker,LimitsImpliesGiesekerGeneral}.

The converse forms the technical core of this paper: we must control a subobject $A_r\hookrightarrow E$ that may vary with $r$ and need not be a coherent subsheaf. The phase expansion for a fixed pair of sheaves gives no uniform threshold for such objects. We first prove that the sequences relevant to destabilization become bounded coherent sequences. The estimates needed for this step are established on $C^n$ and then transferred to the chosen induced family on $X$. Here $C$ is a complex elliptic curve and
\[
 H=\sum_{i=1}^n\operatorname{pr}_i^*[0].
\]
On $C^n$ the central charge has the simpler form
\begin{equation}\label{equation: explicit central charge}
 Z_{r,\theta}(E)=\sum_{j=0}^n\zeta_jr^jb_j(E),\qquad
 \zeta_j=-(-e^{i\theta})^j,\qquad
 b_j(E)=\frac1{j!}\int_{C^n}H^j\ch{n-j}(E).
\end{equation}
Since $\td(C^n)=1$, the $b_j$ are the ordinary Hilbert-polynomial coefficients. In \Cref{Theorem: Dimension determines limit phase}, we show that the HN phases of a bounded family of pure $d$-dimensional sheaves concentrate uniformly at rate $O(r^{-1})$ around the limiting phase
\[
\phi_d = 1 - \frac{d(\pi - \theta)}{\pi}.
\]
This controls every HN phase, not only the phase of the total central charge. For torsion-free sheaves, the two extreme phases can be bounded using presentations and embeddings by line bundles. A pure sheaf of positive codimension has no such embedding into an ambient vector bundle. We instead project its support finitely to projective space of the same dimension, where its direct image is torsion-free, and use a flat degeneration to transfer the resulting phase bounds; see \Cref{proposition: phase bounds for pure sheaves on projective space}.

The next step combines numerical and geometric control. Uniform deformations bound the weighted Hilbert coefficients of potential destabilizers, while a separate mass estimate controls the multiplicities of their ordinary cohomology sheaves. Restriction and induction on support dimension then bound the leading coefficient and bound from below the next coefficient of the coherent image. Grothendieck boundedness applies to its saturation. Finally, a small-mass Hom-vanishing argument forces the original subobject to equal that saturation. This gives the following criterion for sequences in the moving heart.

\begin{theorem}[\Cref{theorem: Narrow inclusions become saturated and coherent}]\label{Introduction: Narrow inclusions become saturated and coherent}
Let $\mathscr B$ be a bounded family of pure $d$-dimensional sheaves, where $1\leq d\leq n$, and let $M>0$ be the constant given by \Cref{Theorem: Dimension determines limit phase}. There is $R>0$ such that, for $E\in\mathscr B$ and $r\geq R$, every exact sequence
    \[
    0\longrightarrow A\longrightarrow E\longrightarrow B\longrightarrow0
    \]
in $\mathscr A_r$ whose nonzero terms have HN phases in $[\phi_d-M/r,\phi_d+M/r]$ is a coherent exact sequence with $A\subset E$ saturated. As $E$, the sequence, and $r\geq R$ vary, the sheaves $A$ and $B$ form bounded families.
\end{theorem}

Uniform comparison of sheaf phases now proves that Gieseker stability implies limit stability on $C^n$. In \Cref{subsec:general-comparison}, we transport potential destabilizing sequences through $\pi^*i_{X*}$ to obtain the comparison on $X$ and complete the proof of \Cref{Thm:MainComparison}.

We also study the opposite sector $0<\theta<\pi/N$, where the limiting heart is perverse coherent. These results suggest further perverse coherent limits in the remaining noncritical directions, and potential stability theory on perverse coherent sheaves; see \Cref{PervCohAsLimits}, \Cref{PervCohAsLimitsGeneral,Conj:LimitOfHearts}. 

Our rays also fit into the framework of quasi-convergent stability conditions developed by Halpern-Leistner--Jiang--Robotis \cite[Definition~2.8]{HLJR2024quasiconvergencestabilityconditions}. Our phase-concentration theorem implies that these rays are quasi-convergent, with associated filtration given by the dimension of the supports of cohomology sheaves \cite[Theorem~2.30]{HLJR2024quasiconvergencestabilityconditions}. The induced stability conditions on the successive quotients are determined by the leading Hilbert-polynomial coefficients, and the rays converge weakly to the corresponding augmented stability condition \cite[Proposition~5.21(1)]{HLR2025spaceaugmentedstabilityconditions}. Our Gieseker comparison captures finer information: although every pure sheaf is limit-semistable in the quasi-convergence sense, eventual semistability along the ray is equivalent to Gieseker semistability.


\subsection{Related works} The relation between Bridgeland stability and classical stability at the large volume limit has been studied from several perspectives. Toda's limit stability on Calabi--Yau threefolds \cite{Toda09LimitStableObjectsOnCY3} relates limiting stability conditions on perverse coherent hearts to curve-counting theories. Relations with counting theories are also investigated, for example in \cite{OPT2020AbelianThreefolds,FeyzbakhshThomas23CurveCountingSDuality,FeyzbakhshThomas23RankRDTfromrank1,FeyzbakhshThomas24RankRDTfromrank0,JLMM2025HigherRankDTPT}.

Bayer–Macrì–Toda identify an object-wise eventual heart in a fixed-$B$ large-volume limit and
prove the corresponding one-way semistability comparison \cite[Lemmas 6.2.1–6.2.2, Proposition 6.2.3]{BMT2014}. For suitable paths in the stability-parameter half-plane for $\mathbb{P}^3$, Pretti identifies coherent sheaves by eventual membership in the moving hearts \cite[Proposition 4.1]{Pretti2021ZeroRank}. 

For surfaces and threefolds, the relation between asymptotic stability and Gieseker stability has been studied from the perspective of wall-crossing and asymptotic geometry of walls, for example in \cite{LoQin2014MiniWalls,Pretti2021ZeroRank,JM22WallsAndAsymptoticsOnThreefolds,JMM23VerticalAsymptotics,bernal2025bridgelandwallsdestabilizingonedimensional}.

There are also some other functorial approaches to Gieseker stability, for example \cite{alvarez-consulFunctorialConstructionModuli2007,heinPostnikovStabilityComplexes2007,InabaModuliOfStableObjects2010}. The use of boundedness and uniform regularity to control relevant subobjects is analogous to Álvarez-Cónsul--King \cite[§5.1, conditions (C:1)–(C:5); §5.2, Definition 5.3 and Theorem 5.10]{alvarez-consulFunctorialConstructionModuli2007}. Their comparison concerns Kronecker modules attached to regular sheaves. \Cref{theorem: Narrow inclusions become saturated and coherent} here instead establishes coherence and saturation for sequences in the moving Bridgeland hearts. Hein–Ploog \cite[Theorem 12]{heinPostnikovStabilityComplexes2007} characterize pure Gieseker semistable sheaves of a fixed Hilbert polynomial by finite homological data in the derived category. Inaba \cite[Definition 3.1 and Example 5.1]{InabaModuliOfStableObjects2010} obtains a related sheaf-moduli construction from strict ample sequences. These use stability structures different from the ray of Bridgeland conditions studied here.

\subsection{Conventions and Notations}\label{paragraph:conventions}
All varieties are over $\mathbb{C}$. We write $n=\dim X$ and reserve $N$ for an ambient projective dimension. On a product $C^n$, $H=\sum_i H_i$ is the product polarization. On a general $X$, the polarization and its embedding in $\mathbb{P}^N$ are fixed throughout.

All functors between derived categories are derived unless stated otherwise.  We write $\mathscr{H}^q(E)$ for standard cohomology sheaves. If $\mathscr{A}$ is a heart, $\mathscr{A}[a,b]$ denotes the extension closure of $\mathscr{A}[j]$ for $a\leq j\leq b$; thus $\Coh(X)[0,n]=\mathrm{D}^{[-n,0]}_{\Coh}(X)$.

The variables $r$ and $\theta$ always denote the radial and angular parameters. When recalling the construction we also use
\begin{equation*}
     \alpha=r\sin\theta,\qquad \beta=r\cos\theta.
\end{equation*}
Once $\theta$ and the variety are fixed, we abbreviate $\sigma_{r,\theta}$, $Z_{r,\theta}$ and $\mathscr{A}_{r,\theta}$ to $\sigma_r$, $Z_r$ and $\mathscr{A}_r$. A phase of a semistable object is denoted by $\phi$; $\phi^\pm$ denote the extreme HN phases. For a nonzero object of a heart, $\varphi=\pi^{-1}\arg Z$ denotes the phase of its total central charge, even when the object is not semistable. These quantities need not agree.

A \textit{bounded family of coherent sheaves} means a family contained, after finite stratification, in a coherent flat family over a finite-type scheme. Single sheaves and finite collections are bounded.

A \textit{bounded family of complexes} in $\Db(X)$ means a collection contained, up to isomorphism, among the derived fibres over $\mathbb{C}$-points of an $S$-perfect complex $\mathcal E$ on $X\times_\mathbb{C} S$ for some finite-type $\mathbb{C}$-scheme $S$, which we may replace by finitely many reduced locally closed strata so that the derived restriction of $\mathcal E$ to each stratum has cohomology sheaves flat over that stratum.

In a statement uniform for a bounded family $\mathscr{B}$, the constants may depend on $\mathscr{B}$, the fixed variety and polarization, the fixed angle $\theta$, and other explicitly fixed data. They are independent of the particular member of $\mathscr{B}$ and of $r$. Uniformity in $\theta$ is not asserted unless stated. Uniformity in a deformation parameter $\lambda$ means uniformity on the specified compact coefficient set. The phrase ``for all sufficiently large $r$'' always refers to a threshold with these dependencies.

\subsection{Structure of the paper} 
Section~\ref{sec:Preliminary} recalls the stability conditions and fixes the Hilbert-polynomial conventions. In Section~\ref{sec:LimitOfHearts}, uniform eventual membership of bounded sheaf families and the converse for fixed objects identify the coherent-heart limit. Section~\ref{sec:LimitImpliesGieseker} proves the implication from limit stability to Gieseker stability. 

For the converse, Section~\ref{sec:uniform-estimates} supplies uniform weighted coefficient estimates and a mass bound for individual coherent cohomology sheaves. Section~\ref{section: Pure sheaf phase estimates and coefficients control} establishes HN-phase concentration for bounded pure families and then tests the coefficients of moving subobjects by deformation. The coherent-image estimate and controlled enlargement yield \Cref{theorem: Narrow inclusions become saturated and coherent}: the relevant heart sequences are bounded saturated coherent sequences. The polynomial comparison then proves the forward implication on bounded families. Finally, Section~\ref{subsec:general-comparison} transfers candidate sequences through $\pi^*i_{X*}$ and applies fixed-class Gieseker boundedness to finish \Cref{Thm:MainComparison}.

The ample-property and opposite-sector results in Section~\ref{sec:LimitOfHearts} are additional consequences of the heart analysis; they are not needed for the stability comparison.

\subsection{AI disclosure}
This project began in March 2026. In mid-June 2026, the author began using OpenAI's GPT-5.5 Pro as a mathematical research assistant. At that stage, the model suggested using Serre's theorem and the closure of torsion classes under quotients in the proof of \Cref{CoherentSheavesInLargeVolumeHeart}. Apart from this contribution, all results and arguments preceding Section 5 were obtained without input from AI models.

The main mathematical difficulty was to establish the implication “Gieseker implies limit”, which required controlling the potentially destabilizing subobjects in varying hearts. Boundedness was recognized early as the central obstacle, and the problem was explored through extensive and iterative discussions with AI models. The key breakthrough came after the release of GPT-6 Astra: with its assistance, a proof was first obtained for torsion-free sheaves and subsequently extended to pure sheaves. The model made substantial contributions to the technical arguments in Sections 5–7, particularly the uniform deformation and mass estimates, the control of Harder--Narasimhan phases, and the coherence and saturation theorem.

The author subsequently substantially simplified, refined, and organized these arguments, and takes full responsibility for the mathematical content and correctness of the manuscript.

\subsection{Acknowledgments}

The author would like to thank Arend Bayer, Yiran Cheng, Hanfei Guo, Alekos Robotis, Benjamin Sung, Yukinobu Toda, Qizheng Yin, and Xiaolei Zhao for enlightening discussions and helpful comments. The author is particularly grateful to Chunyi Li for his uniform support and guidance. Part of this work was carried out during the author’s visit to AMSS and BICMR, whose hospitality is gratefully acknowledged.

The author is supported by the Warwick Mathematics Institute Centre for Doctoral Training, and gratefully acknowledges funding from the University of Warwick. The author was also partially supported by the Royal Society URF\textbackslash R1\textbackslash 201129 and 251025 “Stability condition and application in algebraic geometry”.

\section{Preliminary}\label{sec:Preliminary}
We recall the stability-theoretic conventions, the induced stability conditions on projective varieties, and the polynomial formulation of Gieseker stability. Throughout this section, $X$ is a smooth projective variety of dimension $n$ and $\mathscr{D}$ is a triangulated category. We write $\Db(X)=\Db(\Coh(X))$.

\subsection{Stability conditions}

We use $\phi$ for a phase and reserve $\theta$ for the angle of the ray $z=r e^{i\theta}$. Our conventions for stability conditions follow \cite[Definitions 3.3 and 5.1, Proposition 5.3]{bridgeland_stability_2007} and \cite[Section 2]{Bayer19ShortProofDeformationProperty}.

\begin{definition}\label{Definition:Slicing}
    A \textit{slicing} of $\mathscr{D}$ is a collection of full additive subcategories $\mathscr{P}(\phi)\subset\mathscr{D}$, indexed by $\phi\in\mathbb{R}$, with the following properties:
    \begin{enumerate}[label=(\arabic*)]
        \item\label{Definition:HomologicalShift} $\mathscr{P}(\phi+1)=\mathscr{P}(\phi)[1]$ for every $\phi\in\mathbb{R}$.
        \item If $\phi_1>\phi_2$ and $E_j\in\mathscr{P}(\phi_j)$, then $\Hom_{\mathscr{D}}(E_1,E_2)=0$.
        \item\label{Definition:HNFiltration} Every nonzero object $E\in\mathscr{D}$ admits a finite filtration by distinguished triangles, represented by
        \[
        0=E_0\xrightarrow{s_1}E_1\xrightarrow{s_2}\cdots\xrightarrow{s_m}E_m=E,
        \]
        whose nonzero factors $F_j=\operatorname{Cone}(s_j)$ belong to $\mathscr{P}(\phi_j)$ for real numbers $\phi_1>\cdots>\phi_m$.
    \end{enumerate}
\end{definition}

Fix a free abelian group $\Lambda$ of finite rank and a homomorphism $v:\K_0(\mathscr{D})\twoheadrightarrow\Lambda$. We may replace $\Lambda$ by the image of $v$ when convenient.

\begin{definition}
    A \textit{pre-stability condition} with respect to $(\Lambda,v)$ is a pair $\sigma=(Z,\mathscr{P})$, where $\mathscr{P}$ is a slicing and $Z:\Lambda\to\mathbb{C}$ is a group homomorphism such that
    \[
    Z(v(E))\in\mathbb{R}_{>0}e^{i\pi\phi}
    \qquad\text{for every }0\neq E\in\mathscr{P}(\phi).
    \]
    We call $Z$ the \textit{central charge} and abbreviate $Z(v(E))$ to $Z(E)$.
\end{definition}

\begin{notation}
    A nonzero object of $\mathscr{P}(\phi)$ is \textit{$\sigma$-semistable of phase $\phi$}; we write $\phi_\sigma(E)=\phi$. It is \textit{$\sigma$-stable} if it is simple in the abelian category $\mathscr{P}(\phi)$.
\end{notation}

The filtration in \Cref{Definition:Slicing}\ref{Definition:HNFiltration} is the \textit{Harder--Narasimhan filtration}, or HN filtration, and its factors are uniquely determined up to isomorphism. For a nonzero object $E$ with HN factors $F_1,\ldots,F_m$, set
\[
\phi^+_\sigma(E)=\phi_1,
\qquad
\phi^-_\sigma(E)=\phi_m,
\qquad
m_\sigma(E)=\sum_{j=1}^m|Z(F_j)|.
\]
The last quantity is the \textit{mass} of $E$. When the stability condition is clear, we omit the subscript $\sigma$. For an interval $I\subset\mathbb{R}$, $\mathscr{P}(I)$ denotes the full extension-closed subcategory of objects whose HN factors have phases in $I$, together with the zero object.

\begin{definition}[{\cite[Section 1.2]{kontsevich2008stabilitystructuresmotivicdonaldsonthomas}}]
    A pre-stability condition $\sigma=(Z,\mathscr{P})$ satisfies the \textit{support property} if, for a fixed norm on $\Lambda_{\mathbb{R}}=\Lambda\otimes\mathbb{R}$, there is a constant $C_\sigma>0$ such that
    \[
    |Z(E)|\geq C_\sigma\|v(E)\|
    \qquad\text{for every nonzero $\sigma$-semistable object }E.
    \]
    A pre-stability condition satisfying the support property is called a \textit{stability condition}. This property is independent of the chosen norm.
\end{definition}

We denote the space of such stability conditions by $\Stab_{(\Lambda,v)}(\mathscr{D})$. Its topology is induced by the generalized metric
\begin{equation}\label{definition:distance of stability conditions}
    \begin{aligned}
    \dist(\sigma_1,\sigma_2)
    =\sup_{0\neq E\in\mathscr{D}}\max\biggl\{&
    |\phi^-_{\sigma_1}(E)-\phi^-_{\sigma_2}(E)|,
    |\phi^+_{\sigma_1}(E)-\phi^+_{\sigma_2}(E)|,
    \left|\log\frac{m_{\sigma_2}(E)}{m_{\sigma_1}(E)}\right|
    \biggr\}\in[0,+\infty].
    \end{aligned}
\end{equation}

\begin{theorem}[Deformation theorem, {\cite[Theorem 7.1]{bridgeland_stability_2007}, \cite[Theorem 1.2]{Bayer19ShortProofDeformationProperty}}]\label{Thm:Bridgeland's Deformation Theorem}
    The map
    \[
    \operatorname{Forg}_Z:\Stab_{(\Lambda,v)}(\mathscr{D})\longrightarrow\Hom_{\mathbb{Z}}(\Lambda,\mathbb{C}),
    \qquad (Z,\mathscr{P})\longmapsto Z,
    \]
    is a local homeomorphism. Consequently, when nonempty, $\Stab_{(\Lambda,v)}(\mathscr{D})$ is a complex manifold of dimension $\rk\Lambda$.
\end{theorem}

For $\mathscr{D}=\Db(X)$, we write $\Stab_{(\Lambda,v)}(X)$, or simply $\Stab(X)$ when the lattice is fixed. We also use the equivalent description by a heart and a stability function.

\begin{definition}
    Let $\mathscr{A}$ be an abelian category and $v:\K_0(\mathscr{A})\to\Lambda$ a homomorphism. A \textit{stability function} is a homomorphism $Z:\Lambda\to\mathbb{C}$ such that
    \[
    Z(E)\in\{z\in\mathbb{C}:\im z>0\}\cup\mathbb{R}_{<0}
    \qquad\text{for all }0\neq E\in\mathscr{A}.
    \]
    Define $\varphi(E)=\pi^{-1}\arg Z(E)\in(0,1]$. The object $E$ is \textit{$Z$-semistable} (respectively, \textit{$Z$-stable}) if every nonzero proper subobject $F\subset E$ satisfies $\varphi(F)\leq\varphi(E)$ (respectively, $\varphi(F)<\varphi(E)$).
\end{definition}

A stability function has the \textit{HN property} if every nonzero object admits a finite filtration by subobjects whose nonzero successive quotients are semistable with strictly decreasing phases.

\begin{lemma}[{\cite[Proposition 5.3]{bridgeland_stability_2007}}]
    Giving a pre-stability condition on $\mathscr{D}$ is equivalent to giving a pair $(Z,\mathscr{A})$, where $\mathscr{A}$ is the heart of a bounded $t$-structure on $\mathscr{D}$ and $Z$ is a stability function on $\mathscr{A}$ with the HN property.
\end{lemma}

The correspondence sends $(Z,\mathscr{P})$ to $(Z,\mathscr{P}((0,1]))$. Conversely, the semistable objects in $\mathscr{A}$ determine $\mathscr{P}(\phi)$ for $\phi\in(0,1]$, and the shift rule determines the remaining slices. We use the canonical identification $\K_0(\mathscr{A})\cong\K_0(\mathscr{D})$. The support property can be checked on semistable objects in $\mathscr{A}$, so the same correspondence applies to stability conditions. Thus $\sigma=(Z,\mathscr{A})$ always denotes the description by a heart, whereas $\sigma=(Z,\mathscr{P})$ denotes the description by a slicing.

We use hearts for exact sequences and subobjects, and slicings for HN phases and masses.

\begin{definition}
    A stability condition $\sigma$ on $\Db(X)$ is called \textit{geometric} if for each point $p \in X$, the skyscraper sheaf $\mathcal{O}_p$ is $\sigma$-stable, and all skyscraper sheaves are of the same phase.
\end{definition}

The cohomological amplitude of a geometric stability condition is always bounded, in the following sense.

\begin{lemma}[{\cite[Lemma 10.1]{bridgeland_stability_k3}, \cite[Lemma 2.11]{FLZ2022}}]\label{Lemma:Geometric Amplitude}
    Let $\sigma$ be a geometric stability condition on a smooth projective variety $X$ of dimension $n$, for which skyscraper sheaves are of phase $1$, and let $\mathscr{A} = \mathscr{P}_{\sigma}((0,1])$. Then
    \[
    \mathscr{A} \subset \mathrm{D}^{[-n,0]}_{\Coh} (X),
    \qquad
    \Coh(X) \subset \mathscr{P}_{\sigma}([-n, 1]).
    \]
\end{lemma}





We collect the phase and mass estimates used below.

\begin{lemma}[Phase and mass estimates]\label{lemma:phase-and-mass}
Let $\sigma=(Z,\mathscr P)$ be a stability condition.
\begin{enumerate}[label=(\arabic*), ref=\thelemma~(\arabic*)]
\item\label[lemma]{lemma: Phase Control By Morphism}
For nonzero objects and morphisms $L\twoheadrightarrow F\hookrightarrow U$ in a heart $\mathscr P((a,a+1])$, one has
\[
 \phi^-(L)\leq\phi^-(F)\leq\phi^+(F)\leq\phi^+(U).
\]
For a distinguished triangle $A\to B\to C$ with nonzero terms, one also has
\[
 \phi^-(B)\geq\min\{\phi^-(A),\phi^-(C)\},\qquad
 \phi^+(B)\leq\max\{\phi^+(A),\phi^+(C)\}.
\]
\item\label[lemma]{lemma:charge-phase-mass relation}
If $0\neq B\in\mathscr P([a,b])$ and $\psi\in\mathbb R$, then
\[
\begin{aligned}
 m(B)\min_{t\in[a,b]}\sin\bigl(\pi(t-\psi)\bigr)
 &\leq\im\bigl(e^{-\pi i\psi}Z(B)\bigr)\\
 &\leq m(B)\max_{t\in[a,b]}\sin\bigl(\pi(t-\psi)\bigr).
\end{aligned}
\]
When $[a,b]\subset[\psi-\tfrac12,\psi+\tfrac12]$, the minimum and maximum are attained at $a$ and $b$, respectively.
\item\label[lemma]{lemma:mass-central-charge}
If $0\neq B\in\mathscr P([\psi-\varepsilon,\psi+\varepsilon])$, where $\psi\in\mathbb R$ and $0\leq\varepsilon<\tfrac12$, then
\[
 \cos(\pi\varepsilon)m(B)
 \leq\re\bigl(e^{-\pi i\psi}Z(B)\bigr)
 \leq|Z(B)|\leq m(B).
\]
In particular, mass and the modulus of the central charge are uniformly comparable whenever the HN phases lie in an interval whose width is bounded strictly below $1$.
\end{enumerate}
\end{lemma}
\begin{proof}
The first HN factor of $F$ maps nontrivially to $U$, and the last HN quotient of $F$ is a quotient of $L$. The triangle bounds follow from extension closure of phase intervals. This proves (1). If $B_j$ are the HN factors of $B$, then
\[
 e^{-\pi i\psi}Z(B)
 =\sum_j|Z(B_j)|e^{\pi i(\phi(B_j)-\psi)}.
\]
Taking imaginary parts gives (2); taking real parts and applying the triangle inequality gives (3).
\end{proof}

\subsection{Pullback and pushforward of stability conditions}
We recall the operations used in \cite[Definition 3.1]{li2026remarkstabilityconditionssmooth}. Let $f:Y\to X$ be a morphism of smooth projective varieties. For stability data $\sigma_Y=(Z_Y,\mathscr{P}_Y)$ on $Y$, define candidate induced data on $X$ by
\[
f_\sharp\sigma_Y=(Z_Y\circ f^*,f_\sharp\mathscr{P}_Y),
\qquad
(f_\sharp\mathscr{P}_Y)(\phi) =\{E\in\Db(X):f^*E\in\mathscr{P}_Y(\phi)\}.
\]
Similarly, for stability data $\sigma_X=(Z_X,\mathscr{P}_X)$ on $X$, set
\[
f^\sharp\sigma_X=(Z_X\circ f_*,f^\sharp\mathscr{P}_X),
\qquad
(f^\sharp\mathscr{P}_X)(\phi) =\{E\in\Db(Y):f_*E\in\mathscr{P}_X(\phi)\}.
\]
These formulas do not automatically define pre-stability conditions. Whenever they do, we call them the \textit{pushforward} and \textit{pullback} stability conditions, respectively. All functors here are derived, as in our conventions.

Let $p:C\to C/\langle[-1]\rangle\cong\mathbb{P}^1$ be the degree-two quotient of the elliptic curve. Composing $p^n$ with the symmetric-product quotient gives
\[
\pi:C^n\longrightarrow(\mathbb{P}^1)^n
\longrightarrow\operatorname{Sym}^n(\mathbb{P}^1)\cong\mathbb{P}^n.
\]
This is a finite flat surjective map of degree $2^n n!$, and $\pi^*H_{\mathbb{P}^n}=2H$.

\begin{theorem}[{\cite[Theorem 6.3]{li2026remarkstabilityconditionssmooth},\cite[Theorem 5.6]{LLLMPSZ2026stabilityconditionsmodulispaces}}]\label{StabilityConditionOnPn}
    The stability conditions on $C^n$ in \Cref{StabilityConditionOnEn} induce stability conditions on $\mathbb{P}^n$. For $\alpha\in\mathbb{R}_{>0}$ and $\beta\in\mathbb{R}$, set $z=\beta+i\alpha$ and normalize the induced condition by
    \[
    \tau_z=\bigl((2^nn!)^{-1}(Z_{2z}\circ\pi^*),\pi_\sharp\mathscr{P}_{2z}\bigr).
    \]
    Its central charge is
    \[
    Z_{z,\mathbb{P}^n}(A)=-\int_{\mathbb{P}^n}e^{-zH_{\mathbb{P}^n}}\ch{}(A).
    \]
\end{theorem}
\begin{proof}[Central-charge computation]
    Since $2H=\pi^*H_{\mathbb{P}^n}$, the degree formula and projection formula give
    \[
    Z_{2z}(\pi^*A)
    =-\int_{C^n}\pi^*\bigl(e^{-zH_{\mathbb{P}^n}}\ch{}(A)\bigr)
    =-2^n n!\int_{\mathbb{P}^n}e^{-zH_{\mathbb{P}^n}}\ch{}(A).
    \]
    Multiplication of the central charge by the positive number $(2^n n!)^{-1}$ leaves the slicing unchanged.
\end{proof}

Now fix an embedding $i_X:X\hookrightarrow\mathbb{P}^N$ with $i_X^*H_{\mathbb{P}^N}=H$.

\begin{theorem}[Induced conditions on $X$; {\cite[Theorem 6.5]{li2026remarkstabilityconditionssmooth},\cite[Theorem 6.2]{LLLMPSZ2026stabilityconditionsmodulispaces}}]\label{StabilityConditionOnSmoothProjectiveVarieties}
    For the fixed embedding $i_X$, there is $\alpha_{i_X} > 0$ such that, whenever $\alpha > \alpha_{i_X}$, the induced pair $\sigma_{z,X}=i_X^\sharp\tau_z$ is a stability condition on $X$. Its heart and central charge are
    \[
    \mathscr{A}_{z,X}=\{E\in\Db(X):i_{X*}E\in\mathscr{A}_{z,\mathbb{P}^N}\},
    \]
    \[
    Z_{z,X}(E)=Z_{z,\mathbb{P}^N}(i_{X*}E)
    =-\int_X e^{-zH}\ch{}(E)\frac{\td(X)}{i_X^*\td(\mathbb{P}^N)}.
    \]
\end{theorem}
\begin{proof}[Central-charge computation]
    Grothendieck--Riemann--Roch gives
    \[
    \ch{}(i_{X*}E)\td(\mathbb{P}^N)=i_{X*}\bigl(\ch{}(E)\td(X)\bigr).
    \]
    Substitute this identity into $Z_{z,\mathbb{P}^N}$ and apply the projection formula.
\end{proof}

\subsection{Gieseker stability}
We follow the conventions of \cite[Section 1.2]{HL10GeometryOfModuliOfSheaves}, except that our Hilbert coefficients include the factorials. Throughout $X$ is a smooth projective variety of dimension $n$, $H$ an ample divisor on $X$; in statements involving $i_X: X\hookrightarrow \mathbb{P}^N$, assume that $H = i^*_X H_{\mathbb{P}^N}$ is very ample.

\begin{definition}
    A nonzero coherent sheaf $E$ is \textit{pure} of dimension $d$ on $X$ if any non-zero subsheaf $F\subset E$ has dimension $d$. A pure sheaf of dimension $n = \dim (X)$ will be called \textit{torsion-free}.
\end{definition}

Let $E$ be a coherent sheaf of dimension $d$. Its \textit{Hilbert polynomial} is
\[
P_H(E,m) = \chi(E \otimes \mathcal{O}(mH)) = \sum\limits_{j = 0}^d a_j(E) m ^j,
\qquad 
a_d(E) > 0.
\]
The polynomial has degree $d$; see \cite[Lemma 1.2.1]{HL10GeometryOfModuliOfSheaves}. By Hirzebruch--Riemann--Roch,
\[
a_j(E) = \frac{1}{j!} \int_X H^j(\ch{}(E)\td(X))_{n - j}.
\]
Here the subscript denotes the component of codimension $n - j$. We suppress $H$ when the polarization is fixed and set $a_j(E) = 0$ for $j > d$.
\begin{definition}
    The \textit{reduced Hilbert polynomial} of $E$ is
    \[
    p_H(E,m) \coloneqq \frac{P_H(E,m)}{ a_d(E)}.
    \]
    A pure sheaf $E$ is $H$-Gieseker (semi)stable if every nonzero proper subsheaf $F\subset E$ satisfies
    \[
    p_H(F) <(\leq)  p_H(E).
    \]
    Polynomial inequalities mean inequalities for all sufficiently large $m$, or equivalently lexicographic comparison of coefficients starting with the highest degree.
\end{definition}

Thus Gieseker stability includes purity as part of its definition. All nonzero zero-dimensional sheaves are Gieseker semistable.

On $C^n$, the tangent bundle is trivial, so the Hilbert coefficients agree with the notation fixed in the introduction:
\[
a_j(E) = b_j(E) = \frac{1}{j!} \int_{C^n} H^j\ch{n - j}(E).
\]
We use $a_j$ for the ordinary Hilbert coefficients on a general $X$, $b_j$ for coefficients on $C^n$, and $\tilde{a}_j$ for the coefficients defined next.

\begin{remark}\label{ToddTwistDoesNotAffectGieseker}
    Motivated by \Cref{StabilityConditionOnSmoothProjectiveVarieties}, for the fixed embedding $i_X$, we define the \textit{Todd-twisted Hilbert polynomial} to be 
    \[
    \tilde{P}_H (E,m) \coloneqq \int_X \ch{}(E\otimes\mathcal{O}(mH)) \frac{\td(X)}{i_X^*\td(\mathbb{P}^N)} = \sum\limits_{j = 0}^d \tilde{a}_j(E) m^j . 
    \]
    The reduced Todd-twisted Hilbert polynomial is defined similarly. Since the inverse Todd class of projective space is a polynomial in its hyperplane class, write
    \[
    i^*_X\td(\mathbb{P}^N)^{-1} = \sum\limits_{\ell = 0}^n t_\ell H^\ell,
    \qquad
    t_0 = 1.
    \]
    Comparing coefficients gives the upper-triangular relation
    \begin{equation}
        \tilde{a}_j(E) = \sum\limits_{\ell = 0}^{d - j} t_{\ell} \frac{(j + \ell )!}{j!} a_{j + \ell }(E).
    \end{equation}
    In particular, $\tilde{a}_d(E) = a_d(E)$. If $F$ and $E$ have dimension $d$ and their normalized ordinary coefficients agree above degree $k$, then their normalized twisted coefficients also agree above degree $k$, and
    \[
    \frac{\tilde{a}_k(F)}{\tilde{a}_d(F)} - \frac{\tilde{a}_k(E)}{\tilde{a}_d(E)} = \frac{a_k(F)}{a_d(F)} - \frac{a_k(E)}{a_d(E)}.
    \]
    Hence the sign of the first differing normalized coefficient is unchanged. Reduced Todd-twisted polynomials therefore define exactly Gieseker stability and semistability. This notion would be crucial to the final comparison between limit stability and Gieseker stability in \Cref{subsec:general-comparison}.
\end{remark}

This notation also makes the induced charge transparent:
\[
Z_{r,\theta, X}(E) = - \tilde{P}_H(E,-re^{i\theta}) = \sum\limits_{j = 0}^d \zeta_jr^j \tilde{a}_j(E),
\qquad
\zeta_j = -(-e^{i\theta})^j.
\]

\begin{remark}\label{GiesekerCentralChargeTransformed}
    The central charge expression of Gieseker stability
    \[
    Z_{\text{Gie}}(E; m) \coloneqq -P_H(E,m) + i \tfrac{\partial}{\partial m} P_H(E, m)
    \]
    forms a nice bridge between Gieseker stability and limit stability. For each nonzero $E$, its values lie in the upper half-plane or on the negative real axis for sufficiently large $m$; write \[\phi_{\text{Gie}}(E;m) = \pi^{-1}\arg Z_{\text{Gie}}(E;m).\] If $F$ and $E$ have the same dimension, then
    \[
    p_H(F)\leq p_H(E)
    \quad\Longleftrightarrow\quad
    \phi_{\mathrm{Gie}}(F;m)\leq\phi_{\mathrm{Gie}}(E;m) \quad\text{for }m\gg0,
    \]
    and the same holds with strict inequalities.
    Indeed, suppose that $d=\dim E=\dim F>0$ and that $p_H(F,m)-p_H(E,m)=c m^k+O(m^{k-1})$ with $c\neq0$ and $k<d$. Then
    \[
    \im\!\left(Z_{\mathrm{Gie}}(F;m)\overline{Z_{\mathrm{Gie}}(E;m)}\right)= a_d(F) a_d(E) \bigl((d-k)c\,m^{d+k-1}+O(m^{d+k-2})\bigr).
    \]
    Its eventual sign is the sign of $c$, which is also the sign of the phase difference. If the reduced polynomials agree, the charges are positive multiples of one another; the zero-dimensional case is immediate. Consequently, a pure sheaf is Gieseker (semi)stable exactly when these eventual phase inequalities hold for each nonzero proper subsheaf.

    This formulation fits in the framework of Bayer's polynomial stability conditions \cite{Bayer09PolynomialStabilityCondition}.
\end{remark}

\section{Limit of hearts and the ample property}\label{sec:LimitOfHearts}
We first prove the coherent-heart limit used in the stability comparison. The last two subsections derive the ample property from a horizontal limit and identify the opposite-sector perverse coherent limit. These are additional consequences of the heart analysis. Throughout this section we use the families and parameter convention:
\[
 \beta+i\alpha=re^{i\theta},\qquad
 \mathscr{A}_{r,\theta}=\mathscr{A}_{\alpha,\beta}.
\]
The angle $\theta$ is fixed, while $r$ tends to infinity. The notation
\[
 \lim_{r\to\infty}\mathscr{A}_{r,\theta}=\mathscr{B}
\]
will mean that an object belongs to $\mathscr{B}$ if and only if it belongs to $\mathscr{A}_{r,\theta}$ for all sufficiently large $r$. This is an assertion about eventual membership, rather than convergence in an unspecified topology on $t$-structures.

\subsection{Stability conditions on products}
We recall the construction of stability conditions on products \cite{liu21productstability} needed for the induction on the number of elliptic factors; see also \cite{li2025stabilityconditionsproductscurves,haiden2024stabilitymanifoldetimesetimes}. Let $Y$ be smooth and projective and let $\sigma=(Z,\mathscr{A})$ be a stability condition with respect to $(\Lambda,v)$, with $Z(\Lambda)\subset\mathbb{Q}+i\mathbb{Q}$. The heart is noetherian by \cite[Lemma 12.12]{BLMNPSStabilityConditionInFamilies}. Let $C$ be a smooth projective curve of genus $g$, and choose a degree-one line bundle $\mathcal{O}_C(1)$. Write $p_Y$ and $p_C$ for the projections from $Y\times C$; pushforwards in the following formulas are derived.

The constant family of $t$-structures has heart
\[
\mathscr{A}_C \coloneqq \{ E \in \Db(Y\times C) | p_{Y*}(E\otimes p^*_C\mathcal{O}_C(m)) \in \mathscr{A}  \text{ for } m \gg 0 \};
\]
see \cite[Theorem 2.6.1]{AP06SheavesOftStructures}, also \cite{P07ConstantFamiliesOftStructures}. Following \cite[Section 4.1, equation (4.1)]{li2025stabilityconditionsproductscurves}, define additive functions $a,b,c,d$ by
\begin{equation}
    Z(p_{Y*}(E\otimes p^*_C\mathcal{O}_C(m + g - 1))) = a(E)m + b(E) + i(c(E)m + d(E)).
\end{equation}
The charge factors through
\[
 \widetilde\Lambda_C=\Lambda\oplus(\Lambda/\ker Z),
\]
with class map $\widetilde v_C$ induced by
\begin{align}
    v'_C(E) & \coloneqq v(p_{Y*}(E\otimes p^*_C \mathcal{O}_C(m))) - v(p_{Y*}(E\otimes p^*_C \mathcal{O}_C(m - 1))), \\
    v''_C(E) &\coloneqq v(p_{Y*}(E\otimes p^*_C \mathcal{O}_C(m))) - (m+ 1 - g)\cdot v'_C(E),
\end{align}
for $E\in\Db(Y\times C)$, which are independent of $m\in\mathbb{Z}$.

For $t\in\mathbb{R}_{> 0}$ and $u\in\mathbb{R}$, set 
\begin{equation}\label{Liu'sTiltStability}
    'Z_{t, u, C} (E) \coloneqq a(E) t - d(E) + c(E) u + ic(E)t.
\end{equation}
This is a weak stability function on $\mathscr{A}_C$ with Harder--Narasimhan filtrations \cite[Theorem 3.3 and Lemma 4.3]{liu21productstability}. Its torsion pair is
\begin{align*}
 \mathcal{T}_{t,u,C}
 &=\left\{E\in\mathscr{A}_C\ \middle|\
 \begin{array}{l}
 \text{every nonzero quotient }E\twoheadrightarrow F\text{ satisfies}\\
 c(F)=0\text{ or }\re{}'Z_{t,u,C}(F)<0
 \end{array}\right\},\\
 \mathcal{F}_{t,u,C}
 &=\left\{E\in\mathscr{A}_C\ \middle|\
 \begin{array}{l}
 \text{every nonzero subobject }F\hookrightarrow E\text{ satisfies}\\
 c(F)>0\text{ and }\re{}'Z_{t,u,C}(F)\geq0
 \end{array}\right\}.
\end{align*}
Define
\[
 \mathscr{A}_{t,u,C}=\langle\mathcal{T}_{t,u,C},\mathcal{F}_{t,u,C}[1]\rangle,
 \qquad
 Z_{s,t,u,C}(E)=b(E)+s c(E)-u a(E)
                  +i\bigl(d(E)-t a(E)-u c(E)\bigr).
\]

\begin{theorem}[{\cite[Theorem 5.9]{liu21productstability};\cite[Theorem 4.1]{li2025stabilityconditionsproductscurves}}]\label{LiuStabilityOnProduct}
For $s, t\in\mathbb{R}_{>0}$ and $u\in\mathbb{R}$, the pair $\sigma_{s, t, u, C} \coloneqq (Z_{s, t, u, C},\mathscr{A}_{t, u, C})$ is a stability condition on $\Db(Y\times C)$ with respect to $(\tilde{\Lambda}_C,\tilde{v}_C)$.   
\end{theorem}

We now specialize to $C^n=C_1\times\cdots\times C_n$, with each $C_i$ elliptic, $H_i = \text{pr}^*_i[0]$ and $H=H_1+\cdots+H_n$. The lattice we would be using is
\[
\Lambda_n = (\mathbb{Z} \oplus \mathbb{Z})^{\otimes n},
\qquad
v_I(T) = \int_{C^n} \left(\prod_{i\in I}H_i \right)\ch{n - |I|}(T);
\]
see \cite[Section~4.2 and Theorem~4.5]{li2025stabilityconditionsproductscurves}. The real-parameter family is defined over this same lattice by \cite[Proposition~4.5]{LLLMPSZ2026stabilityconditionsmodulispaces}. Each coefficient $b_j$ is a fixed linear form on $\Lambda_{n,\mathbb R}$.

\begin{proposition}[{\cite[Lemma 4.4 and Theorem 4.5]{li2025stabilityconditionsproductscurves},\cite[Proposition 4.5]{LLLMPSZ2026stabilityconditionsmodulispaces}}]\label{StabilityConditionOnEn}
For $(\beta,\alpha)\in\mathbb{R}\times\mathbb{R}_{>0}$, there is a stability condition on $C^n$ 
\[
 \sigma_{\alpha,\beta,n}=(Z_{\alpha,\beta,n},\mathscr{A}_{\alpha,\beta,n}),
 \qquad
 Z_{\alpha,\beta,n}(E)=-\int_{C^n}e^{-(\beta+i\alpha)H}\ch{}(E).
\]
Its heart is constructed inductively from $\mathscr{A}_{\alpha,\beta,1}=\Coh(C_1)$. If $p_{C^{n - 1}}:C^n\to C^{n-1}$ is the projection and $p_C:C^n\to C$ is the last projection, set
\[
 \mathscr{C}_{\alpha, \beta, n} \coloneqq \{ A\in \Db(C^{n - 1}\times  C) | p_{C^{n - 1}*}(A\otimes p^*_C\mathcal{O}_C(m))\in \mathscr{A}_{\alpha, \beta, n-1}  \text{ for all } m\gg 0 \}. 
\]
Applying \Cref{LiuStabilityOnProduct} with $s=t=\alpha$ and $u=\beta$ gives a torsion pair $(\mathcal{T}_{\alpha,\beta,n},\mathcal{F}_{\alpha,\beta,n})$ in this heart and
\[
 \mathscr{A}_{\alpha,\beta,n}
 =\langle\mathcal{T}_{\alpha,\beta,n},\mathcal{F}_{\alpha,\beta,n}[1]\rangle.
\]
Skyscraper sheaves are stable of phase $1$.
\end{proposition}

\subsection{Coherent sheaves as the limiting heart}
Fix $\theta\in (\pi -\pi/n,\pi)$. We prove below
\[
\lim_{\substack{
r\to\infty
}}
\mathscr{A}_{r,\theta,n} =  \Coh(C^n).
\]
The statements then are generalized to all smooth projective varieties. We first prove uniform eventual membership for bounded families of coherent sheaves, and then the converse for an arbitrary fixed derived object. For the next lemma to be precise, we set $\mathcal{T}_{\alpha,\beta,1} = \Coh(C)$.

\begin{lemma}\label{StructureSheafInHeart}
    For any $\theta\in(\pi - \tfrac{\pi}{n},\pi)$, $\mathcal{O}_{C^n} \in \mathcal{T}_{\alpha,\beta, n}$, and hence $\mathcal{O}_{C^n} \in \mathscr{A}_{\alpha,\beta,n}$.
\end{lemma}
\begin{proof}
    We prove by induction on $n$. The case $n = 1$ is trivial. Assume $\mathcal{O}_{C^{n - 1}} \in \mathscr{A}_{\alpha,\beta,n - 1}$. 
    For $m > 0$,
    \[
    p_{C^{n - 1}*}(\mathcal{O}_{C^{ n}} \otimes p^*_C\mathcal{O}_{C}(m)) = \mathcal{O}_{C^{n - 1}}^{\oplus m},
    \]
    so the induction hypothesis puts $\mathcal{O}_{C^n}$ in $\mathscr{C}_{\alpha,\beta, n}$. The intermediate heart satisfies
    \[
    \mathscr{C}_{\alpha,\beta,n}\subset \langle\mathscr{A}_{\alpha,\beta,n}[-1],\mathscr{A}_{\alpha,\beta,n}\rangle.
    \]
    By \cite[Corollary 2.16]{FLZ2022}, $\mathcal{O}_{C^n}$ is $\sigma_{\alpha, \beta, n}$-stable, so its phase lies in $(-1,1]$. Since $H^n = n!$,
    \begin{equation}
        Z_{r,\theta, n}(\mathcal{O}_{C^n}) = - \frac{H^n}{n!}(-re^{i\theta})^n = r^n e^{\pi i - n(\pi - \theta)i}.
    \end{equation}
    The phase in $(-1,1]$ is consequently $1-n(\pi-\theta)/\pi\in(0,1)$. Thus the structure sheaf belongs to both hearts, whose intersection is $\mathcal{T}_{\alpha,\beta,n}$.
\end{proof}

Tensoring by $\mathcal{O}_{C^n}(kH)$ carries $\mathscr{A}_{\alpha,\beta,n}$ to $\mathscr{A}_{\alpha,\beta+k,n}$. In particular, for each fixed integer $m$, the sheaf $\mathcal{O}_{C^n}(-mH)$ belongs to $\mathcal{T}_{\alpha,\beta,n}$ for all sufficiently large $r$: apply \Cref{StructureSheafInHeart} to $(\alpha,\beta+m)$, whose argument tends to $\theta$.

In the following, we fix the projection $p = p_{C^{n - 1}}: C^n = C^{n - 1}\times C \rightarrow C^{n - 1}$ into the first $n - 1$ factors, and let $i: C^{n - 1} \hookrightarrow C^n$ be the embedding given by $x \mapsto (x,0)$.

\begin{proposition}\label{CoherentSheavesInLargeVolumeHeart}
    Fix a bounded family of coherent sheaves $\mathscr{B}$ on $C^n$ and $\theta\in(\pi - \pi/n,\pi)$. There is a constant $R = R_{\mathscr{B},\theta}$ such that
    \[
    A\in\mathscr{A}_{r,\theta,n} 
    \qquad
    (A\in \mathscr{B}, r\geq R_{\mathscr{B},\theta}).
    \]
    The threshold can also be chosen uniformly for $\theta$ in a fixed compact subinterval of $(\pi-\pi/n,\pi)$. At rational construction parameters and for $n\geq2$, every such $A$ lies in $\mathcal{T}_{\alpha,\beta,n}$.
\end{proposition}
\begin{proof}
    We prove by induction on $n$. The case $n = 1$ is immediate since $\mathscr{A}_{\alpha,\beta, 1} = \Coh(C)$. Suppose it holds in dimension $n-1$, uniformly for angles in compact subintervals of the stated sector. Fix such a compact interval $J\subset(\pi-\pi/n,\pi)$; it also lies in the corresponding sector in dimension $n-1$. We first prove uniform membership in the intermediate heart $\mathscr{C}_{\alpha,\beta,n}$, then use a uniform presentation of sheaves and the fact that the torsion class is closed under quotients to conclude the induction.
    
    Relative Serre vanishing for the bounded family gives one integer $M$ such that
    \[
    A_m \coloneqq p_* (A \otimes p^*_C\mathcal{O}_{C}(m))
    \]
    is a coherent sheaf for every $A\in\mathscr{B}$ and $m\geq M$. The sheaves $A_M$, $\mathscr{H}^{-1}(i^*A)$, and $\mathscr{H}^{0}(i^*A)$ form bounded families. The induction hypothesis therefore gives one radius after which $A_M$ lies in $\mathscr{A}_{r,\theta,n-1}$ and
    \[
    i^*A\in
    \mathrm{D}^{[-1,0]}_{\mathscr{A}_{r,\theta,n-1}}.
    \]
    Here and below, superscripts on $\mathrm{D}_{\mathscr{A}}$ indicate cohomological degrees for the heart $\mathscr{A}$; in particular, an object of $\mathscr{A}$ shifted by $[1]$ has $\mathscr{A}$-cohomology in degree $-1$.

    We have triangles from the divisor sequence
    \[
    A_m\longrightarrow A_{m+1}\longrightarrow i^*A\longrightarrow A_m[1].
    \]
    If $A_m\in\mathscr{A}_{r,\theta,n-1}$, this triangle puts $A_{m+1}$ in $\mathrm{D}^{[-1,0]}_{\mathscr{A}_{r,\theta,n-1}}$. On the other hand, the inclusion $\mathscr{A}_{r,\theta,n-1}\subset\mathrm{D}^{\leq0}_{\Coh}$ implies $\Coh\subset\mathrm{D}^{\geq0}_{\mathscr{A}_{r,\theta,n-1}}$. Since $A_{m+1}$ is a sheaf, it belongs to the heart. We have shown that:
    \[
    A_M\in\mathscr{A}_{r,\theta,n - 1} \implies A_m\in\mathscr{A}_{r,\theta,n - 1} \text{ for all } m\geq M.
    \]
    Hence $A\in\mathscr{C}_{\alpha,\beta,n}$, uniformly in $A$.
    
    Next choose a single integer $m_0$ such that every $A(m_0H)$ in the family is globally generated, with higher cohomology zero; this is possible since a bounded family has a uniform regularity bound. The evaluation sequences are
    \[
    0\longrightarrow K_A\longrightarrow H^0(A(m_0H))\otimes\mathcal{O}_{C^n}(-m_0H) \longrightarrow A\longrightarrow0.
    \]
    The kernels $K_A$ form a bounded family: after a finite stratification of a parameter space for $\mathscr{B}$, these are kernels of the relative evaluation maps. Applying the intermediate-heart argument to this bounded family of kernels, we may assume $K_A \in \mathscr{C}_{\alpha,\beta,n}$. By \Cref{StructureSheafInHeart}, we might also enlarge the radius $r$ so that $\mathcal{O}_{C^n}(-m_0H) \in \mathcal{T}_{\alpha,\beta,n}$. This last choice is uniform for $\theta\in J$, since the arguments of $(\beta+m_0)+i\alpha$ stay in the required sector for all sufficiently large $r$. 
    
    All three terms then lie in $\mathscr{C}_{\alpha,\beta,n}$, so the displayed sequence is also exact in that heart. Since the torsion class is closed under quotients, we have $A\in\mathcal{T}_{\alpha,\beta,n}$.
\end{proof}

\begin{remark}\label{Remark:HorizontalHeartLimit}
    The same induction applies with $\alpha = \alpha_0 > 0$ fixed and $\beta\rightarrow-\infty$: every bounded family of sheaves belongs to $\mathscr{A}_{\alpha_0,\beta,n}$ once $-\beta$ exceeds a family-dependent threshold. 
\end{remark}

\begin{proposition}\label{LimitImpliesCoherent}
Fix $\theta\in(\pi-\pi/n,\pi)$. If $E\in\Db(C^n)$ belongs to $\mathscr{A}_{r,\theta,n}$ for all sufficiently large $r$, then $E$ is a coherent sheaf.
\end{proposition}
\begin{proof}
By \Cref{Lemma:Geometric Amplitude}, such an object has no positive standard cohomology. If its least nonzero cohomology degree is $-m<0$, standard truncation gives a nonzero morphism
\[
 \mathscr{H}^{-m}(E)[m]\longrightarrow E
\]
which induces the identity on $\mathscr{H}^{-m}(E)$. Apply \Cref{CoherentSheavesInLargeVolumeHeart} to the singleton family $\{\mathscr{H}^{-m}(E)\}$. At a sufficiently large radius, both $\mathscr{H}^{-m}(E)$ and $E$ lie in the same heart, contradicting
\[
\Hom(\mathscr{A}_{r,\theta,n}[m],\mathscr{A}_{r,\theta,n})=0
\qquad(m>0). \qedhere
\]
\end{proof}

\begin{corollary}\label{Cor:UniformHeartOnX}
Let $i_X:X\hookrightarrow\mathbb{P}^N$ be the chosen embedding. Fix $\theta\in(\pi-\pi/N,\pi)$. For every bounded family $\mathscr{B}\subset\Coh(X)$ there is $R_{\mathscr{B},\theta}$ such that
\[
 \mathscr{B}\subset\mathscr{A}_{r,\theta,X}
 \qquad(r\geq R_{\mathscr{B},\theta}).
\]
Moreover, an object of $\Db(X)$ belongs to $\mathscr{A}_{r,\theta,X}$ for all sufficiently large $r$ if and only if it is a coherent sheaf. The same assertions hold on $\mathbb{P}^N$.
\end{corollary}
\begin{proof}
The functor $\pi^*i_{X*}$ preserves bounded families of sheaves and detects standard cohomology, because $i_{X*}$ is exact and faithful and $\pi$ is faithfully flat. The descriptions of the induced hearts in \Cref{StabilityConditionOnPn,StabilityConditionOnSmoothProjectiveVarieties} therefore reduce the claims to the two preceding propositions on $C^N$, with the positive rescaling of the radius dictated by $\pi^*H_{\mathbb{P}^N}=2H$.
\end{proof}

\subsection{Ample property}\label{subsection: ample property}

Tensoring by $\mathcal{O}_X(H)$ shifts $\beta$ by $1$. Thus the horizontal limit in \Cref{Remark:HorizontalHeartLimit} has the following formulation. 

\begin{definition}
A stability condition $\sigma=(Z,\mathscr{A})$ on a polarized variety $(X,H)$ has the \emph{ample property} if every coherent sheaf $E$ satisfies $E(mH)\in\mathscr{A}$ for all sufficiently large integers $m$. It has the \emph{uniform ample property} if, for every bounded family $\mathscr{B}\subset\Coh(X)$, there is one integer $m_{\mathscr{B}}$ such that this holds for every $E\in\mathscr{B}$ and every $m\geq m_{\mathscr{B}}$.
\end{definition}

\begin{remark}
For a heart with the ample property, eventual membership of $E(mH)$ also forces $E$ to be a sheaf. Indeed, all skyscraper sheaves belong to the heart because twisting does not change them. Every object of the heart consequently has no positive standard cohomology: a highest nonzero cohomology sheaf in degree $j>0$ would give a nonzero map to $\mathcal{O}_x[-j]$, contradicting the heart axioms. If $E$ has negative cohomology, apply the ample property to its lowest nonzero cohomology sheaf and use the truncation morphism as in \Cref{LimitImpliesCoherent}.
\end{remark}

\begin{corollary}
All stability conditions constructed in \cite{li2026remarkstabilityconditionssmooth} have the uniform ample property.
\end{corollary}

\begin{corollary}
For each of these hearts there is an equivalence
\[
 \Db(\mathscr{A}_{\alpha,\beta,X})\simeq\Db(X).
\]
More generally, this holds for a bounded heart in $\Db(X)$ with the ample property.
\end{corollary}
\begin{proof}
Use the realization functor obtained from the standard filtered derived category as in \cite[Appendix A.2, A.5--A.6, pp. 38--39]{Beilinson1987DerivedPerverse}. Its triangulated construction in this enhancement is justified by \cite[Theorem 3.11, Remark 3.10 and Appendix Theorem A.3]{PsaroudakisVitoria2018Realisation}. It is enough to prove essential surjectivity, by \cite[Theorem B, equivalently Theorem 2.9]{CHZ19DerivedEquivalences}.

Represent an arbitrary object by a bounded complex $E^\bullet$ of coherent sheaves. Choose $q>0$ such that $L=\mathcal{O}_X(qH)$ is globally generated and $\mathscr{H}^i(E)(jqH)\in\mathscr{A}$ for every term $\mathscr{H}^i(E)$ and every $j\geq1$. There are only finitely many terms, so the ample property supplies one such $q$. Put $V=H^0(X,L)$ and $\ell=\dim V$. The evaluation map gives the exact dual Koszul complex
\begin{equation}\label{KoszulComplex}
 0\longrightarrow\mathcal{O}_X\longrightarrow V^\vee\otimes L
 \longrightarrow\bigwedge^2V^\vee\otimes L^2
 \longrightarrow\cdots\longrightarrow\bigwedge^\ell V^\vee\otimes L^\ell
 \longrightarrow0.
\end{equation}
It is locally split, since the generating sections have no common zero. Tensoring with $E^\bullet$ therefore preserves exactness, and the resulting total complex gives
\[
 E^\bullet\simeq
 \operatorname{Tot}\left[E^\bullet\otimes V^\vee\otimes L
 \longrightarrow E^\bullet\otimes\bigwedge^2V^\vee\otimes L^2
 \longrightarrow\cdots\right],
\]
with the first Koszul term in degree $0$. Every term belongs to $\mathscr{A}\cap\Coh(X)$. Viewed as a bounded complex in $\mathscr{A}$, its realization is this total complex; this compatibility follows from the construction using the filtration by complex degree. Thus every object lies in the essential image.
\end{proof}

The ample property puts skyscraper sheaves in the heart but does not require their stability. For geometric stability conditions, we ask whether the converse principle holds.

\begin{conjecture}
Every geometric stability condition on a smooth polarized variety $(X,H)$, normalized so that skyscraper sheaves have phase $1$, has the ample property with respect to $H$.
\end{conjecture}

\subsection{Perverse coherent sheaves as the limiting heart}

We use perversities depending only on support dimension, following \cite[Section 3.1]{Bayer09PolynomialStabilityCondition}. For more details on the theory of perverse coherent sheaves, we refer the reader to \cite{AB10PerverseCoherentSheaves,Ka04TStructureOnDerivedCategories}.

\begin{definition}
    A function $p: \{0,1,\cdots,n\} \rightarrow \mathbb{Z}$ is called a \textit{perversity function} if $p$ is monotone decreasing, and if $\overline{p}: \{0,1,\cdots,n\} \rightarrow \mathbb{Z}$ (called the \textit{dual perversity}) given by $\overline{p}(d) = -d - p(d)$ is also monotone decreasing.
\end{definition}

In other words we require that $p(d) \geq p(d + 1) \geq p(d) - 1$. Given a perversity function in the above sense, the function $X^{\text{top}} \rightarrow \mathbb{Z}$, $x\mapsto p(\dim x)$ is a monotone and comonotone perversity function in the sense of \cite{AB10PerverseCoherentSheaves}. We set $\dim:X^{\text{top}} \rightarrow \mathbb{Z}, x\mapsto \dim(x)$.

Let $\mathcal{A}^{p,\leq k}$ be the following increasing filtration of $\Coh (X)$ by abelian subcategories:
\[
\mathcal{A}^{p,\leq k} = \left\{ F \in \Coh(X) \mid p(\dim \text{supp}  F ) \geq -k  \right\}.
\]

\begin{theorem}[\cite{AB10PerverseCoherentSheaves,Ka04TStructureOnDerivedCategories}]
 If $p$ is a perversity function, then the following pair defines a bounded $t$-structure on $\Db (X)$:
\[
\mathcal D^{p,\leq 0} = \left\{ E\in\Db(X) \mid \mathscr{H}^{-k}(E) \in \mathcal{A}^{p,\leq k} \text{ for all } k\in\mathbb{Z}   \right\},
\]
\[
\mathcal D^{p,\geq 0} = \left\{ E\in\Db(X) \mid \Hom (A, E) = 0 \text{ for all } k\in\mathbb{Z} \text{ and } A \in \mathcal{A}^{p,\leq k}[k + 1]  \right\}.
\]
\end{theorem}

Objects in the heart $\mathcal{A}^p = \mathcal D^{p,\leq 0} \cap \mathcal D^{p,\geq 0}$ are called \textit{perverse coherent sheaves}, and we denote it by $\Coh^p(X)$.

Recall that on a smooth projective variety $X$ of dimension $n$, we have the Grothendieck dual functor:

\begin{equation}
    \mathbb{D}_X: \Db(X) \rightarrow \Db(X)^{\text{op}}, A \mapsto \RHom(A, \omega_X[n]).
\end{equation}

The dual $\mathbb{D}_X$ yields an anti-equivalence between $\Coh^p(X)$ and $\Coh^{\overline{p}}(X)$. In particular, 
\[
\Coh^0(X) = \Coh(X), 
\qquad
\mathbb{D}_X(\Coh(X)) = \Coh^{-\dim}(X),
\]
where $\Coh^{-\dim}(X)$ is the category of Cohen-Macaulay complexes on $X$; see \cite[Section 6]{YeZhang06DualizingComplexesAndPerverseSheaves}. This fact will be used in the proof of \Cref{PervCohAsLimits}.

For a stability condition $\sigma = (Z,\mathscr{P})$, duality gives $\mathbb{D}_X(\sigma) = (W,\mathscr{Q})$, where
\[
 W(E)=\overline{Z(\mathbb{D}_X(E))},\qquad
 \mathscr{Q}(\phi)=\mathbb{D}_X\bigl(\mathscr{P}(-\phi)\bigr).
\]

\begin{remark}
    Without the complex conjugation on the central charge and minus sign on the slicing, the stability condition obtained would be a stability condition on $\Db(X)^{\text{op}}$, rather than $\Db(X)$.
\end{remark}

On $C^n$, it is convenient to shift the resulting slicing by two so that skyscraper sheaves again have phase $1$.
\begin{lemma}\label{DualStabilityConditions}
    For every real phase $\phi$,
    \[
    \overline{Z_{r,\theta,n}(\mathbb{D}_{C^n}E)} = Z_{r,\pi-\theta, n},
    \qquad
    \mathbb{D}_{C^n}(\mathscr{P}_{r,\theta,n}(\phi)) = \mathscr{P}_{r,\pi - \theta, n}(2 - \phi).
    \]
\end{lemma}
\begin{proof}
    For $v = \sum\limits_{j = 0}^n v_j\in H^*(C^n,\mathbb{Q})$, $v_j\in H^{2j}(C^n,\mathbb{Q})$, define $P(v) \coloneqq \sum\limits_{j = 0}^n (-1)^j v_j$. Since $\omega_{C^n} \cong \mathcal{O}_{C^n}$,
    \[
    P(e^{mH}) = e^{-mH},
    \qquad
    \ch{}(\mathbb{D}_{C^n}A) = (-1)^n P(\ch{}(A)), 
    \qquad
    \int_{C^n}P(v) = (-1)^n \int_{C^n}v.
    \]
    Direct computation shows that
    \begin{align}
        \overline{Z_{\alpha, \beta, n}(\mathbb{D}_{C^n}(A))} &= \overline{-\int_{C^n} e^{-(\beta + i\alpha)H} \cdot \ch{}(\mathbb{D}_{C^n}A)} \\
        & = -\int_{C^n} e^{-(\beta - i\alpha)H} \cdot \ch{}(\mathbb{D}_{C^n}(A)) \\
        & = -\int_{C^n} e^{-(\beta - i\alpha)H} \cdot (-1)^nP(\ch{}(A)) \\
        & = -(-1)^n\int_{C^n} P(e^{(\beta - i \alpha)H}\ch{}(A)) \\ 
        & = - \int_{C^n}e^{-(-\beta + i\alpha)H} \ch{}(A)  = Z_{\alpha, -\beta, n}(A)
    \end{align}
    Since $\mathbb{D}_{C^n}$ is a contravariant equivalence, it sends semistable objects to semistable objects. The uniqueness statement \cite[Theorem 1.1]{li2025stabilityconditionsproductscurves}, combined with the central charge computation and skyscraper phase one convention, implies that there is an integer $k$ such that
    \[
    \mathbb{D}_{C^n}(\mathscr{P}_{\alpha,\beta, n} (\phi)) = \mathscr{P}_{\alpha,-\beta, n} (2k - \phi)
    \]
    Since duality fixes skyscraper sheaves, we conclude that $k = 1$.
\end{proof}

We will also need to control the endpoint of the phase interval. This step cannot be replaced by formally interchanging a limit and a tilt.

\begin{lemma}\label{Lemma:PhaseOneEndpoint}
Let $\sigma=(Z,\mathscr{P})$ be geometric, with skyscraper phase $1$. Every $T\in\mathscr{P}(1)$ has $\mathscr{H}^0(T)$ supported in dimension zero. Consequently, if $F$ is a coherent sheaf with no nonzero zero-dimensional subsheaf and $F\in\mathscr{P}(0,1]$, then $F\in\mathscr{P}(0,1)$.
\end{lemma}
\begin{proof}
The category $\mathscr{P}(1)$ has finite length. Each closed point in the support of $\mathscr{H}^0(T)$ gives a nonzero morphism $T\to\mathcal{O}_x$: use $T\in\mathrm{D}^{\leq0}_{\Coh}$ and the quotient $\mathscr{H}^0(T)\twoheadrightarrow\mathcal{O}_x$. Since $\mathcal{O}_x$ is simple in $\mathscr{P}(1)$, it is a Jordan--H\"older factor of $T$. There can be only finitely many such points.

For any such $T$ and any sheaf $F$ as in the statement, every morphism $T\to F$ factors through $\mathscr{H}^0(T)$, by standard truncation. It is therefore zero. A phase-one first Harder--Narasimhan factor of $F$ is impossible.
\end{proof}

Define
\[
 q(0)=0,\quad q(d)=1-d\ (d>0),
 \qquad
 \overline q(0)=0,\quad\overline q(d)=-1\ (d>0).
\]
Also put
\[
 \mathcal{T}_0=\Coh_0(C^n),
 \qquad
 \mathcal{F}_0=\{F\in\Coh(C^n):\Hom(\Coh_0(C^n),F)=0\}.
\]
Then $\Coh^{\overline q}(C^n)=\langle\mathcal{F}_0[1],\mathcal{T}_0\rangle$.

\begin{proposition}\label{PervCohAsLimits}
    For fixed $\theta\in(0,\tfrac{\pi}{n})$, we have 
    \[
    \lim\limits_{r\rightarrow\infty} \mathscr{A}_{r, \theta, n} = \Coh^{q}(C^n).
    \]
    The forward membership assertion is uniform on bounded families of complexes in $\Coh^q(C^n)$.
\end{proposition}
\begin{proof}
Write $\vartheta=\pi-\theta\in(\pi-\pi/n,\pi)$ and
$\mathscr{B}_r=\mathscr{P}_{r,\vartheta,n}[1,2)$. We prove its eventual-membership limit directly.

By \Cref{CoherentSheavesInLargeVolumeHeart,Lemma:PhaseOneEndpoint}, a bounded family in $\mathcal{F}_0$ belongs uniformly to $\mathscr{P}_{r,\vartheta,n}(0,1)$ for sufficiently large $r$. Every zero-dimensional sheaf belongs to $\mathscr{P}_{r,\vartheta,n}(1)$ for every $r$. Thus a bounded family in $\langle\mathcal{F}_0[1],\mathcal{T}_0\rangle$ belongs uniformly to $\mathscr{B}_r$; apply the preceding statements to its bounded families of standard cohomology sheaves.

Conversely, suppose a fixed $E$ belongs to $\mathscr{B}_r$ for all large $r$. Since
\[
 \mathscr{B}_r=
 \langle\mathscr{P}_{r,\vartheta,n}(0,1)[1],\mathscr{P}_{r,\vartheta,n}(1)\rangle,
\]
the geometric amplitude bound and \Cref{Lemma:PhaseOneEndpoint} imply that $E\in\mathrm{D}^{\leq0}_{\Coh}$ and $\mathscr{H}^0(E)$ is zero-dimensional. If its lowest nonzero degree is $-m$ with $m\geq2$, choose $r$ large enough that $\mathscr{H}^{-m}(E)$ belongs to $\mathscr{A}_{r,\vartheta,n}$. The nonzero truncation morphism
$\mathscr{H}^{-m}(E)[m]\to E$ contradicts phase orthogonality: its source has all phases greater than $m\geq2$, whereas $E$ has all phases less than $2$.
Finally, a nonzero zero-dimensional subsheaf $T\subset\mathscr{H}^{-1}(E)$ would give a nonzero morphism $T[1]\to E$. Its source has phase $2$, another contradiction. Hence $E\in\langle\mathcal{F}_0[1],\mathcal{T}_0\rangle$.

We have proved $\lim\mathscr{B}_r=\Coh^{\overline q}(C^n)$. By \Cref{DualStabilityConditions}, $\mathbb{D}_{C^n}(\mathscr{B}_r)=\mathscr{A}_{r,\theta,n}$. Duality identifies the limiting heart with $\Coh^q(C^n)$ and preserves bounded families of complexes on the smooth projective variety.
\end{proof}


Let $i_X: X \hookrightarrow \mathbb{P}^N$ be a projective embedding. Note that $N$ is not necessarily the dimension of $X$.

\begin{lemma}
Let $p:\{0,\ldots,N\}\to\mathbb Z$ be a support-dimension perversity. On $X$ we use its restriction to dimensions at most $\dim X$. For the chosen embedding $i_X:X\hookrightarrow\mathbb{P}^N$, we have
\[
 E\in\Coh^p(X)\quad\Longleftrightarrow\quad
 \pi^*i_{X*}E\in\Coh^p(C^N).
\]
\end{lemma}
\begin{proof}
Closed pushforward and finite faithfully flat pullback preserve and detect the dimensions of supports of standard cohomology sheaves, so they preserve and detect $\mathrm{D}^{p,\leq0}$. Apply the same statement to duality and the dual perversity to obtain the assertion for $\mathrm{D}^{p,\geq0}$. Here duality commutes with closed pushforward; its compatibility with finite flat pullback involves tensoring by the invertible relative dualizing sheaf, which does not change support conditions.
\end{proof}

Combined with \Cref{Cor:UniformHeartOnX}, we conclude that 

\begin{proposition}\label{PervCohAsLimitsGeneral}
    For  fixed $\theta\in(0,\tfrac{\pi}{N})$, we have 
    \[
    \lim\limits_{r\rightarrow\infty} \mathscr{A}_{r, \theta, X} = \Coh^{q}(X),
    \qquad
    \lim\limits_{r\rightarrow\infty} \mathscr{A}_{r, \pi - \theta, X} = \Coh(X).
    \]
\end{proposition}

The remaining directions suggest a family of support-dimension perversities. One must exclude the rays on which a leading term lies on the boundary of the chosen phase interval.

\begin{conjecture}\label{Conj:LimitOfHearts}
Fix $\theta\in(0,\pi)$ such that $d(\pi-\theta)/\pi\notin\mathbb{Z}$ for every $1\leq d\leq n$, and set
\[
 p_\theta(d)=-\left\lfloor\frac{d(\pi-\theta)}{\pi}\right\rfloor.
\]
Then
\[
 \lim_{r\to\infty}\mathscr{A}_{r,\theta,n}=\Coh^{p_\theta}(C^n).
\]
\end{conjecture}

\begin{remark}
The exclusion of critical rays matters. For example, on $C^2$ at $\theta=\pi/2$ the product heart is the usual slope tilt at $\beta=0$. Thus $\mathcal{O}_{C^2}(mH)$ with $m>0$ stays in the heart without a shift, whereas the formula $p_{\pi/2}(2)=-1$ gives a perverse heart containing these line bundles shifted by $[1]$. On these boundary rays, lower coefficients of the central charge enter the limiting heart.
\end{remark}

\begin{figure}[htbp]
\centering
\begin{tikzpicture}[scale=1.1]

\def\n{4}
\pgfmathsetmacro{\ang}{180/\n}
\pgfmathsetmacro{\leftang}{180-\ang}
\def\R{4}

\fill[blue!18] 
  (0,0) -- (0:\R)
  arc[start angle=0,end angle=\ang,radius=\R]
  -- cycle;

\fill[red!18]
  (0,0) -- (\leftang:\R)
  arc[start angle=\leftang,end angle=180,radius=\R]
  -- cycle;

\draw[dashed, thick] (0,0) -- (\ang:\R);
\draw[dashed, thick] (0,0) -- (\leftang:\R);

\draw[->, thick] (-4.4,0) -- (4.4,0) node[right] {$\beta$};
\draw[->, thick] (0,0) -- (0,3.4) node[above] {$\alpha$};

\node[blue!60!black] at (1.9,0.45) {$\operatorname{Coh}^{q}(C^n)$};
\node[red!60!black] at (-1.9,0.45) {$\operatorname{Coh}(C^n)$};

\node at (2.65,0.95) {$0<\arg<\frac{\pi}{n}$};
\node at (-2.55,0.95) {$\pi-\frac{\pi}{n}<\arg<\pi$};



\end{tikzpicture}
\caption{Schematic limiting heart sectors with $n = 4$}
\end{figure}
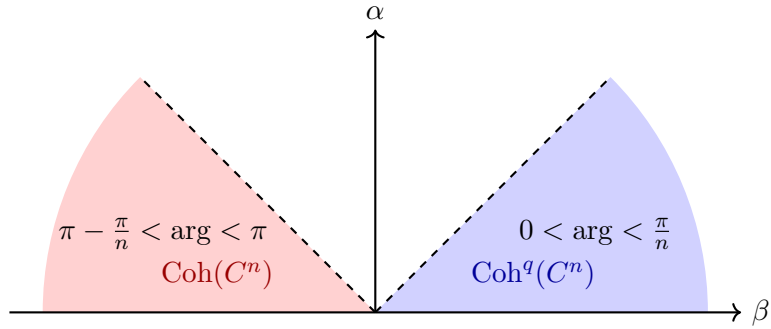

\section{Limit implies Gieseker}\label{sec:LimitImpliesGieseker}

We now compare stability along a fixed ray with Gieseker stability. Throughout the discussion on $C^n$, fix
\[
    \theta\in(\pi-\pi/n,\pi),
    \qquad \delta\coloneqq\pi-\theta\in(0,\pi/n),
\]
and let $r\to\infty$. With the normalization in
\eqref{equation: explicit central charge},
\begin{equation}\label{eq:limit-gieseker-ray-charge}
    Z_{r,\theta}(F)=\sum_{j=0}^n\zeta_j r^j b_j(F),
    \qquad
    \zeta_j=-e^{-ij\delta}=e^{i(\pi-j\delta)},
    \qquad
    b_j(F)=\frac{1}{j!}\int_{C^n}H^j\ch{n-j}(F).
\end{equation}
Since $\td(C^n)=1$, these are precisely the Hilbert-polynomial coefficients:
$P_H(F,m)=\sum_j b_j(F)m^j$.

The phase comparison needed below is a comparison of polynomials. We record it separately to make the distinction between pointwise and uniform assertions explicit.

\begin{lemma}[Asymptotic comparison of sheaf phases]\label{lemma:asymptotic-sheaf-phase-comparison}
    Let $\mathscr{B}_1$ and $\mathscr{B}_2$ be bounded families of nonzero coherent sheaves on $C^n$. There is $R>0$ such that every member of either family belongs to $\mathscr{A}_{r,\theta}$ for $r\geq R$. For such a sheaf $F$, write
    \[
        \varphi_{r,\theta}(F)\coloneqq\frac{1}{\pi}\arg Z_{r,\theta}(F)\in(0,1].
    \]
    This is the phase of its central charge; no semistability is assumed. After increasing $R$, the following hold for every $F\in\mathscr{B}_1$, $E\in\mathscr{B}_2$, and $r\geq R$:
    \begin{enumerate}
        \item If $\dim F<\dim E$, then $\varphi_{r,\theta}(F)>\varphi_{r,\theta}(E)$.
        \item If $\dim F=\dim E$, then the comparisons $<,=,>$ between $\varphi_{r,\theta}(F)$ and $\varphi_{r,\theta}(E)$ agree with those between $p_H(F)$ and $p_H(E)$.
    \end{enumerate}
\end{lemma}

\begin{proof}
    Uniform membership in the heart follows from \Cref{CoherentSheavesInLargeVolumeHeart}. If $F$ has dimension $d$, then $b_j(F)=0$ for $j>d$ and $b_d(F)>0$, so
    \[
        Z_{r,\theta}(F)=\zeta_d b_d(F)r^d+O(r^{d-1}),
        \qquad
        \varphi_{r,\theta}(F)=1-\frac{d\delta}{\pi}+O(r^{-1}).
    \]
    For $d=0$, the phase is exactly $1$. This proves the first assertion.

    For the second assertion, suppose both sheaves have dimension $d$ and set
    \[
        u_j(F)=\frac{b_j(F)}{b_d(F)},
        \qquad
        u_j(E)=\frac{b_j(E)}{b_d(E)}.
    \]
    If $p_H(F)\ne p_H(E)$, let $k<d$ be the largest index for which $u_k(F)\ne u_k(E)$. Then
    \begin{align*}
        \frac{\operatorname{Im}\bigl(Z_{r,\theta}(F)\overline{Z_{r,\theta}(E)}\bigr)}{b_d(F)b_d(E)}
        &=(u_k(F)-u_k(E))\operatorname{Im}(\zeta_k\overline{\zeta_d})r^{d+k}
          +O(r^{d+k-1})\\
        &=(u_k(F)-u_k(E))\sin((d-k)\delta)r^{d+k}
          +O(r^{d+k-1}).
    \end{align*}
    The sine is strictly positive because $0<(d-k)\delta<\pi$. Since both charges lie in the semi-closed upper half-plane, the sign of this imaginary part is the sign of their phase difference. It is therefore the sign of the first nonzero coefficient of $p_H(F)-p_H(E)$, read in descending degree.

    If $p_H(F)=p_H(E)$, then
    \[
        \frac{Z_{r,\theta}(F)}{b_d(F)}
        =\frac{Z_{r,\theta}(E)}{b_d(E)}
    \]
    for every $r$, and their phases are equal whenever both sheaves lie in the heart. This also covers $d=0$.

    Finally, a bounded family has only finitely many Hilbert polynomials. The possible charge polynomials and pairwise comparisons above thus form finite sets, so their thresholds can be chosen uniformly over $\mathscr{B}_1$ and $\mathscr{B}_2$.
\end{proof}

\begin{proposition}\label{LimitsImpliesGieseker}
    Fix $\theta\in(\pi-\pi/n,\pi)$, and let $0\ne E\in\Db(C^n)$. Suppose that there is $R_E>0$ such that $E\in\mathscr{A}_{r,\theta}$ is $\sigma_{r,\theta}$-semistable for every $r\geq R_E$. Then $E$ is a pure $H$-Gieseker semistable sheaf. If $E$ is $\sigma_{r,\theta}$-stable for every $r\geq R_E$, then it is $H$-Gieseker stable.
\end{proposition}

\begin{proof}
    By \Cref{LimitImpliesCoherent}, $E$ is a coherent sheaf. Let $d=\dim E$. For each fixed nonzero proper subsheaf $F\subset E$, put $Q=E/F$. The finite family $\{F,E,Q\}$ is bounded, so \Cref{CoherentSheavesInLargeVolumeHeart} gives a threshold after which all three sheaves lie in $\mathscr{A}_{r,\theta}$. Their coherent short exact sequence is then also a short exact sequence in this heart.

    If $E$ is not pure, we can choose $F$ with $\dim F<d$. By \Cref{lemma:asymptotic-sheaf-phase-comparison}, this subobject would have strictly larger phase than $E$ for all sufficiently large $r$, contrary to semistability. Thus $E$ is pure.

    If $E$ is not Gieseker semistable, there would be a nonzero proper subsheaf $F\subset E$ with $p_H(F)>p_H(E)$. Purity gives $\dim F=d$, and \Cref{lemma:asymptotic-sheaf-phase-comparison} again makes $F$ a destabilizing subobject for all sufficiently large $r$. This proves Gieseker semistability.

    Finally, suppose that $E$ is eventually Bridgeland stable. If it is not Gieseker stable, the semistability just proved would give a nonzero proper subsheaf $F\subset E$ with $p_H(F)=p_H(E)$. The normalized central charges are then identical. Once the coherent sequence lies in $\mathscr{A}_{r,\theta}$, the proper subobject $F$ has the same phase as $E$, contradicting stability.
\end{proof}

\begin{remark}
    The proof of \Cref{LimitsImpliesGieseker} only needs to place one chosen subsheaf and its quotient in the heart at a time. The collection of all subsheaves of a bounded family need not itself be bounded, and no threshold uniform over that collection is asserted.
\end{remark}

For a general smooth projective variety, the same argument works using the Todd-twisted Hilbert polynomial from \Cref{ToddTwistDoesNotAffectGieseker}.

\begin{corollary}\label{LimitsImpliesGiesekerGeneral}
    Let $X$ be a smooth projective variety with an embedding $i:X\hookrightarrow\mathbb{P}^N$ and $H=i^*\mathcal{O}_{\mathbb{P}^N}(1)$, and fix $\theta\in(\pi-\pi/N,\pi)$. Suppose that the stability conditions $\sigma_{r,\theta,X}$ of \Cref{StabilityConditionOnSmoothProjectiveVarieties} are defined for all sufficiently large $r$ along this ray. If $0\ne E\in\Db(X)$ belongs to $\mathscr{A}_{r,\theta,X}$ and is $\sigma_{r,\theta,X}$-semistable for all sufficiently large $r$, then $E$ is a pure $H$-Gieseker semistable sheaf. Eventual stability likewise implies $H$-Gieseker stability.
\end{corollary}

The reverse implication requires a different argument. A subobject in the varying heart can be a complex rather than a subsheaf, so Gieseker stability alone does not directly exclude every possible Bridgeland destabilizer. The remaining sections provide the uniform estimates needed for this step on bounded families.

\section{Uniform deformations and control of mass and support}\label{sec:uniform-estimates}
Throughout the remaining sections, we work on $C^n$ with fixed $\theta\in(\pi-\pi/n,\pi)$. We use the coefficients $b_j$ and $\zeta_j$ of \eqref{equation: explicit central charge}, so $|\zeta_j|=1$ and $\zeta_0=-1$. The $b_j(T)$ are integers. Until the last subsection, we work on $X = C^n$.

We first deform the coefficient tuple $\zeta=(\zeta_0,\ldots,\zeta_n)$ while keeping the constant coefficient equal to $-1$. The estimates below use only the scalar coefficients $b_j$. The support property itself is always taken with respect to the lattice $\Lambda_n$; we do not replace that lattice by the tuple $(b_0,\ldots,b_n)$. With the notions in \Cref{StabilityConditionOnEn}, we have
\[
b_j(T) = \sum\limits_{|I| = j} v_I(T).
\]

The weights $r^j$ in the support estimate are essential: an unweighted estimate on a compact interval would not by itself give the required uniformity as $r\to\infty$.

\begin{lemma}\label{Lemma:Deformations with uniform support constant}
    For $\zeta = (\zeta_0,\cdots,\zeta_n)$ in \eqref{equation: explicit central charge}, and fixed $\varepsilon > 0$, there is a compact neighborhood $V$ of $\zeta$ in $\{\lambda_0 = -1\}$ and a continuous family $\{\sigma_{r,\lambda}\}$ for $r\geq 1$ and $\lambda\in V$ and a constant $C>0$ such that
    \[
    Z_{r,\lambda}(T)=\sum_{j=0}^n\lambda_jr^jb_j(T),
    \qquad
    \sigma_{r,\zeta}=\sigma_r,
    \qquad
    \dist(\sigma_{r,\lambda},\sigma_r)<\varepsilon.
    \]
    Here $\dist$ is the metric in \eqref{definition:distance of stability conditions}, including the mass term. Every skyscraper sheaf has phase $1$ in this family. If $T$ is $\sigma_{r,\lambda}$-semistable, then
    \begin{equation}\label{equation: Support property 1}
    \sum_{j=0}^n r^j|b_j(T)|\leq C|Z_{r,\lambda}(T)|.
    \end{equation}
    Consequently, for every $T\in\Db(X)$,
    \begin{equation}\label{equation: Support property on coefficients}
    \sum_{j=0}^n r^j|b_j(T)|\leq Cm_{r,\lambda}(T),
    \qquad 
    |b_j(T)|\leq Cr^{-j}m_{r,\lambda}(T).
    \end{equation}
    After shrinking $V$, the arguments in $(0,\pi)$ define continuous functions
    \[
    \psi_j(\lambda)=\pi^{-1}\arg\lambda_j,
    \qquad
    1>\psi_1(\lambda)>\cdots>\psi_n(\lambda)>0.
    \]
\end{lemma}
\begin{proof}
    The product family is continuous for all real parameters by \cite[Proposition~4.5]{LLLMPSZ2026stabilityconditionsmodulispaces}. Apply the deformation theorem to the compact subfamily $\{\sigma_t:1\leq t\leq4\}$, with deformed charges $\sum_j\lambda_jt^jb_j$. All numerical stability conditions on $X$ are geometric by \cite[Corollary~2.15]{FLZ2022}; hence the condition $\lambda_0=-1$, together with continuity, keeps the skyscraper phase equal to $1$. Uniqueness from the central charge and this phase, \cite[Theorem~1.1]{li2025stabilityconditionsproductscurves}, identifies the local lifts on overlaps. A finite collection of deformation charts therefore gives a common coefficient neighborhood, which we shrink to a compact neighborhood $V$. We may arrange that the metric distance to $\sigma_t$ is smaller than $\varepsilon$ throughout $[1,4]\times V$.

The local support-property theorem \cite[Theorem~1.2]{Bayer19ShortProofDeformationProperty}, followed by a finite compactness argument, gives a common support constant on this family for a fixed norm on the construction lattice. Since each $b_j$ is a linear form on that lattice and $1\leq t\leq4$, it follows that
\[
 \sum_{j=0}^n t^j|b_j(T)|\leq C|Z_{t,\lambda}(T)|
 \quad\text{for every $\sigma_{t,\lambda}$-semistable $T$}.
\]

We extend this estimate to every radius by multiplication isogenies. Let $f_m=[m]:X\to X$. By Grothendieck--Riemann--Roch, we have
\[
 b_j(f_{m*}T)=m^{2j}b_j(T),\qquad
 Z_{t,\lambda}(f_{m*}T)=Z_{m^2t,\lambda}(T).
\]
By \cite[Lemma~5.3]{li2026remarkstabilityconditionssmooth}, $f_m^\sharp\sigma_{t,\lambda}$ are stability conditions for every deformed charge under consideration. Their skyscraper phase is $1$. We define $\sigma_{m^2t,\lambda}\coloneqq f_m^\sharp\sigma_{t,\lambda} $. The same uniqueness theorem shows that this definition agrees for different representations $r=m^2t$ and agrees with the original family at $\lambda=\zeta$. The intervals $[m^2,4m^2]$ cover $[1,\infty)$, so these compatible local definitions yield a continuous family for every $r\geq1$.

By the definition of $f_m^\sharp$, an object $T$ is $\sigma_{m^2t,\lambda}$-semistable if and only if $f_{m*}T$ is $\sigma_{t,\lambda}$-semistable. Applying the exact functor $f_{m*}$ to HN filtrations gives, for every nonzero $T$,
\begin{equation}\label{eq:IsogenyPhasesAndMass}
 \phi^\pm_{m^2t,\lambda}(T)=\phi^\pm_{t,\lambda}(f_{m*}T),
 \qquad m_{m^2t,\lambda}(T)=m_{t,\lambda}(f_{m*}T).
\end{equation}
Thus induction does not increase the metric distance. It also transfers the preceding scalar support estimate, since $r^jb_j(T)=t^jb_j(f_{m*}T)$. This proves \eqref{equation: Support property 1} and the uniform metric bound. Summing over HN factors proves \eqref{equation: Support property on coefficients}. Finally, the strict inequalities between the arguments hold at $\zeta$, hence on a sufficiently small $V$.
\end{proof}

We retain a slightly larger coefficient neighborhood when constructing the family and then take $V$ compactly contained in its interior. When $r\gg 0$, twisting by a fixed line bundle changes the slicing very mildly.

\begin{lemma}\label{lemma: fixed twists}
For every finite set $K\subset\mathbb Z$, there are constants $M, R>0$ such that, for $r\geq R$, $\lambda\in V$, $k\in K$, and every nonzero $T\in\Db(X)$,
\[
|\phi^\pm_{r,\lambda}(T(kH))-\phi^\pm_{r,\lambda}(T)|<M/r.
\]
\end{lemma}
\begin{proof}
    The charge $Z_{r,\lambda}(T(-kH))$ equals $Z_{r,\lambda'}(T)$, where
    \[
    \lambda'_h = \sum^{h}_{j = 0} \binom{h}{j}\lambda_{j} (-k/r)^{h - j}.
    \]
    For $r\geq 1$ these coefficients satisfy $|\lambda'_h -\lambda_h| \leq C_k/r$, with $C_k$ uniform in the compact coefficient set. The assertion then follows from Bridgeland's deformation theorem.
\end{proof}

The coefficient estimate \eqref{equation: Support property on coefficients} controls the numerical class of a complex, in which cohomology classes occur with alternating signs. We now bound the positive sum of the top multiplicities of its individual standard cohomology sheaves. For example, $G\oplus G[1]$ has zero numerical class, although both nonzero cohomology sheaves have the same support as $G$. The next estimate has no bounded-family hypothesis and can therefore be applied to moving subobjects and restriction errors.

\begin{theorem}\label{theorem:Individual cohomological multiplicities}
    There is $A > 0$ such that, for every $r\geq 1$, $\lambda\in V$, $0\leq d \leq n$, and $T\in \Db(X)$ with $\dim \mathrm{Supp} \mathscr{H}^q (T) \leq d$ for all $q$,
    \begin{equation}
        d! r^d \sum_q b_d(\mathscr{H}^q (T)) \le A m_{r,\lambda} (T).
    \end{equation}
    The constant is independent of $d$, $T$, and the standard cohomological amplitude of $T$.
\end{theorem}

\begin{proof}
    First take $\lambda = \zeta$, and we omit the angle subscript when it is clear from the context. Choose an integer $a\geq 4$ such that $aH$ is very ample, and define
    \[
    C =\max_{\substack{0\leq j\leq n \\ 1\leq t \leq 4 }} \exp \dist(\sigma_{te^{i\theta}+ ja}, \sigma_{te^{i \theta}}).
    \]
    The maximum is finite since the indicated family of stability conditions is compact. By twist covariance and the definition of the metric, we have
    \begin{equation}
        m_{\sigma_{te^{i\theta}}} (F(-jaH)) = m_{\sigma_{te^{i\theta} + ja}}(F) \leq C m_{\sigma_{te^{i\theta}}}(F) \quad (0\leq j \leq n, 1\leq t\leq 4).
    \end{equation}
    Write $r = m^2 t$, where $m = \lfloor \sqrt{r}\rfloor$ and $1\leq t < 4$, and set $F = f_{m*}T$. The isogeny identities and exactness of finite pushforward yield
    \begin{equation}\label{eq:numerical relation in Individual cohomological multiplicities}
        m_{\sigma_{t}}(F) = m_{r}(T), 
        \qquad
        b_d(\mathscr{H}^q(F)) = m^{2d}b_d(\mathscr{H}^q (T)).
    \end{equation}
    Choose divisors $D_1,\ldots,D_d\in |aH|$ successively so that, for each $k=1,\ldots,d$, a local defining equation of $D_k$ is a non-zero-divisor both on
    \[
    \mathcal O_{D_1\cap\cdots\cap D_{k-1}} \qquad
    \text{ and } \qquad
    \mathscr{H}^q(F)\big|_{D_1\cap\cdots\cap D_{k-1}}
    \]
    for every $q$, where the intersection for $k=1$ is understood to be $X$. Such choices are possible because $|aH|$ is base-point-free and, at each step, there are only finitely many associated points of these restricted sheaves to avoid. Set
    \[
        D \coloneqq D_1\cap\cdots\cap D_d,
    \]
    and let $F|_{D}$ be the derived restriction of $F$ to $D$. By the regularity condition, restriction commutes with standard cohomology:
    \[
        \mathscr{H}^q(F|_{D})
        \cong \mathscr{H}^q(F)\big|_{D}
        \qquad\text{for every }q.
    \]
    Since $\dim\operatorname{Supp}\mathscr{H}^q(F)\le d$ for every $q$, each $\mathscr{H}^q(F|_{D})$ is either zero or zero-dimensional. Every zero-dimensional sheaf is semistable of phase $1$ and has mass equal to its length, so the nonzero factors $\mathscr{H}^q(F|_D)[-q]$ of the standard Postnikov filtration have strictly decreasing phases $1-q$. Consequently,
    \[
    m_{\sigma_{t}}(F|_D) = \sum_q \mathrm{length}\text{ } \mathscr{H}^q (F|_D) = d! a^d \sum_q b_d(\mathscr{H}^q (F)).
    \]
    On the other hand, the Koszul resolution and mass subadditivity \cite[Proposition 3.3]{Ike21MassGrowth} give
    \[
    m_{\sigma_{t}}(F|_D) \leq \sum_{j = 0}^d \binom{d}{j} m_{\sigma_{t}}(F(-jaH)) \leq 2^dC m_{\sigma_{t}}(F).
    \]
    Combining all these and \eqref{eq:numerical relation in Individual cohomological multiplicities}, we have
    \[
    d! r^d \sum_q b_d(\mathscr{H}^q (T)) \leq  (\frac{2t}{a})^d C m_r (T) \leq 2^n Cm_r(T).
    \]
    For general $\lambda\in V$, the metric bound in \Cref{Lemma:Deformations with uniform support constant} implies $m_r(T)\leq e^{\varepsilon}m_{r,\lambda}(T)$. Thus $A\coloneqq e^{\varepsilon}2^nC$ works for the whole family.
\end{proof}

The following point-wise form is useful when the object itself varies with $r$.

\begin{corollary}\label{Corollary: Small mass forces support drop}
Let $A$ be as in \Cref{theorem:Individual cohomological multiplicities}, and fix $M>0$ and $0\leq k\leq n$. For every $r>\max\{1,AM\}$, every $\lambda\in V$, and every $T\in\Db(X)$ with $m_{r,\lambda}(T)\leq Mr^k$, one has
    \[
    \dim\operatorname{Supp}\mathscr{H}^q(T)\leq k\quad\text{for every }q,
    \qquad
    \sum_q b_k(\mathscr{H}^q(T))\leq\frac{AM}{k!}.
    \]
    No relation between the objects considered at different radii is required.
\end{corollary}
\begin{proof}
The assertion is immediate for $T=0$. Otherwise, let $d$ be the largest dimension of the support of a cohomology sheaf of $T$. Its total degree $d!\sum_qb_d(\mathscr{H}^q(T))$ is a positive integer. The theorem therefore implies $r^d\leq AMr^k$. Since $r>\max\{1,AM\}$, this forces $d\leq k$. Applying the theorem with $d=k$ gives the second inequality.
\end{proof}

\section{Pure sheaf phase estimates and coefficient control}\label{section: Pure sheaf phase estimates and coefficients control}

Keep $\theta\in(\pi-\pi/n,\pi)$ fixed, and shrink the neighborhood in \Cref{Lemma:Deformations with uniform support constant} to a compact polydisc $V$ about $\zeta$. Write
    \[
    \phi_d=\pi^{-1}\arg\zeta_d=1-\frac{d(\pi-\theta)}{\pi},
    \qquad
    \psi_d(\lambda)=\pi^{-1}\arg\lambda_d.
    \]
Thus $\psi_d(\zeta)=\phi_d$. All error terms and thresholds in statements about sheaves are uniform over the specified bounded family and over $\lambda\in V$. $\phi_{r,\lambda}$ and $\phi^\pm_{r,\lambda}$ retain their semistable and HN meanings.

We first show that every HN phase of a pure $d$-dimensional sheaf approaches $\psi_d(\lambda)$. The leading term of its central charge determines only the total-charge phase, so this requires a separate argument. We then bound the coefficients of potential destabilizing subobjects.

\subsection{Pure sheaf phase estimates}

\begin{theorem}\label{Theorem: Dimension determines limit phase}
    For every bounded family $\mathscr{B}$ of pure $d$-dimensional sheaves on $C^n$, $1\leq d \leq n$, there are uniform constants $M = M_{\mathscr{B}}, R_{\mathscr{B}}, C_{\mathscr{B}}$ such that
    \[
    |\phi^{\pm}_{r,\lambda}(F) - \psi_d(\lambda) | \leq M/r, 
    \qquad 
    m_{r,\lambda}(F) \leq C_{\mathscr{B}}\cdot r ^d,
    \qquad
    (F\in\mathscr{B}, r\geq R_{\mathscr{B}}, \lambda\in V).
    \]
    In particular $|\phi^{\pm}_r(F) - \phi_d|  \leq  M/r$.
\end{theorem}

We first prove the estimates for line bundles supported on $Y = C^d \hookrightarrow C^n$, then for torsion-free sheaves supported on a linear subspace $\mathbb{P}^d \hookrightarrow \mathbb{P}^n$. Using a generization argument, we then prove the estimates for pure sheaves on $\mathbb{P}^n$, after which \Cref{Theorem: Dimension determines limit phase} follows.

\begin{remark}\label{remark: Phase comparison for torsion-free sheaves}
    For a torsion-free sheaf $E$, the estimate follows from \Cref{lemma: Phase Control By Morphism} and the morphisms
\[
 \mathcal{O}(-mH)^{\oplus N'} \twoheadrightarrow E \hookrightarrow \mathcal{O}(mH)^{\oplus N'}
\]
for $m,N'$ very big, which eventually become an epimorphism and a monomorphism also in the heart $\mathscr{A}_r$. The embedding comes from dualizing a presentation of $E^\vee$ and using $E\hookrightarrow E^{\vee\vee}$ as $E$ is torsion-free. 

The line-bundle sandwich applies to torsion-free sheaves, but a pure sheaf of positive codimension cannot embed into an ambient vector bundle. We therefore pass through projective space: a finite projection to the support dimension makes the direct image torsion-free, and a flat degeneration realizes this sheaf on a linear subspace. Phase bounds on the special fiber pass to the projectively equivalent fibers at nonzero parameters. The adjunction maps for $\pi$ then return these bounds to $C^n$.
\end{remark}

\begin{remark}\label{remark on phases of pushforward and pullback}
    As recalled in \Cref{StabilityConditionOnPn}, the stability conditions on $C^n$ descend to  stability conditions $\tau_{r,\lambda} = \pi_{\sharp} \sigma_{r,\lambda} $ on $\mathbb{P}^n$ through $\pi:C^n\rightarrow \mathbb{P}^n$, and
    \begin{equation}
        \phi^{\pm}_{\tau_{r,\lambda}}(F) = \phi^{\pm}_{r,\lambda}(\pi^* F).
    \end{equation} 
    These conditions are invariant under $\operatorname{PGL}_{n+1}$ by \cite[Corollary 3.5.2]{P07ConstantFamiliesOftStructures}. On the other hand, pushforward does not necessarily preserve stability. Instead, for $A\in\Db(C^n)$, we have
    \begin{equation}
        \phi^-_{\sigma}\left(A\otimes \mathcal{O}((-2n - 2)H)\right)=\phi_{\sigma}^-(A\otimes \omega^{-1}_\pi) \leq \phi^-_\tau(\pi_*A) \leq \phi^+_{\tau}(\pi_* A) \leq \phi^+_\sigma (A).
    \end{equation}
    Indeed, for the first HN factor $B_+$ and the last HN factor $B_-$ of $\pi_*A$, by adjunction we have
    \begin{align*}
        0 &\neq \Hom(B_+,\pi_*A) \cong \Hom(\pi^*B_+, A), \\
        0 &\neq \Hom(\pi_*A,B_-) \cong \Hom(A\otimes \omega^{-1}_\pi, \pi^*B_- ).
    \end{align*}
\end{remark}

\begin{lemma}\label{CoordinateLineBundles}
Let $j:Y=C^d\hookrightarrow C^n$ be a coordinate inclusion. For a bounded family of line bundles $L$ on $Y$, the sheaves $j_*L$ are $\sigma_{r,\lambda}$-stable and satisfy
\begin{equation}\label{CoordinateLineBundlesEstimate}
 \phi_{r,\lambda}(j_*L)=\psi_d(\lambda)+O(r^{-1})
\end{equation}
uniformly over the family and over $\lambda\in V$.
\end{lemma}
\begin{proof}
Fourier--Mukai transforms on the point factors send $j_*L$ to a line bundle on the full product. Every line bundle is a simple semihomogeneous bundle, hence is stable for every numerical stability condition by \cite[Corollary 2.16]{FLZ2022}. Apply this assertion to the transformed stability condition and then invert the equivalence.

The support has dimension $d$, and $b_d(j_*L)=1$. Thus $b_j(j_*L)=0$ for $j>d$ and
    \[
    Z_{r,\lambda}(j_*L)=\lambda_dr^d+O(r^{d-1}).
    \]
The lower coefficients take only finitely many values in a bounded family, so the error is uniform. At $\lambda=\zeta$, \Cref{CoherentSheavesInLargeVolumeHeart} selects the phase lift in $(0,1]$. Stability, continuity, and the fixed positive distance of $\psi_d(V)$ from $0$ and $1$ select the same lift throughout a sufficiently small $V$.
\end{proof}

\begin{lemma}\label{lemma:Support on linear subspace}
    Let $i :\mathbb{P}^d \hookrightarrow \mathbb{P}^n$ be a linear inclusion. For a bounded family $\mathscr{B}$ of nonzero torsion-free sheaves $V$ on $\mathbb{P}^d$, one has uniformly
    \begin{equation}
        \phi^{\pm}_{r,\lambda} (\pi^*i_* V) = \psi_d (\lambda) + O(r^{- 1}).
    \end{equation}
    The same assertion holds for $i_*V$ in the descended slicing.
\end{lemma}
\begin{proof}
    We first prove the assertion for $G = \mathcal{O}_{\mathbb{P}^d}(k)$ with $k$ fixed. Fix a coordinate inclusion $j: Y = C^d\times \{0\}^{n - d} \hookrightarrow X = C^n$; we might pick the linear subspace $\mathbb{P}^d$ such that the following diagram commutes:
    \[
    \begin{tikzcd}[row sep=large, column sep=large]
    Y
    \arrow[r, hook, "j"]
    \arrow[d, "\pi_d"']
    &
    X
    \arrow[d, "\pi = \pi_n"]
    \\
    \operatorname{Sym}^{d}\mathbb{P}^{1} = \mathbb{P}^{d}
    \arrow[r, hook, "i"]
    &
    \mathbb{P}^{n}
    \end{tikzcd}
    \]
    where $\pi_d$ is the usual double-symmetric quotient of degree $2^d d!$. Define
    \[
    E_k = i_* G,
    \qquad
    A_k = j_*(\pi_d)^* G = j_*\mathcal{O}_Y (2kH_Y),
    \]
    where $H_Y = H|_Y$. Note that $\mathcal{O}_{\mathbb{P}^d}$ is a direct summand of $\pi_{d*}\mathcal{O}_Y$, hence $E_k$ is a direct summand of $\pi_{*}A_k$.
    Therefore
    \[
    \phi^-_{\sigma}(\pi^*\pi_{*}A_k) \leq \phi^-_{\sigma} (\pi^* E_k) \leq \phi^+_{\sigma} (\pi^* E_k) \leq \phi^+_{\sigma}(\pi^*\pi_{*}A_k).
    \]
    It remains to bound the phase of $\pi_{*}A_k$. By \Cref{remark on phases of pushforward and pullback}, we have
    \[
    \phi_{\sigma}(A_{k -n - 1}) \leq \phi^-_{\tau}(\pi_* A_k)= \phi^-_{\sigma}(\pi^*\pi_{*}A_k)\leq \phi^+_{\sigma}(\pi^*\pi_{*}A_k)=\phi^+_{\tau}(\pi_* A_k) \leq \phi_{\sigma}(A_k),
    \]
    and the assertion follows from \Cref{CoordinateLineBundles}.

    For a torsion-free sheaf $V$ on $\mathbb{P}^d$, choose a presentation and an embedding as in \Cref{remark: Phase comparison for torsion-free sheaves}
    \begin{equation}\label{eq:TorsionFreeSheavesOnLinearSubspace}
        \mathcal{O}_{\mathbb{P}^d}(-a)^{\oplus N'} \twoheadrightarrow V \hookrightarrow \mathcal{O}_{\mathbb{P}^d}(b)^{\oplus N'}.
    \end{equation}
    Since $\pi^* i_*$ is exact on coherent sheaves, the bound for line bundles on $\mathbb{P}^d$ then bounds the phases of $\pi^* i_* V$ by \Cref{lemma: Phase Control By Morphism}.

    For bounded $V$, the integers and multiplicities can be chosen uniformly, hence the error $O(r^{-1})$ can also be made uniform for $\mathscr{B}$.
\end{proof}

Using the open-heart property and generic flatness as in \cite{AP06SheavesOftStructures}, we then transfer phase bounds from a special fiber to the general fibers. Its rationality hypothesis is used only to obtain noetherian rotated hearts.

\begin{lemma}\label{Lemma:Phase estimate from special fiber to general fiber}
    Fix a numerical stability condition $\tau = (Z,\mathscr{P}) $ on $\mathbb{P}^n$ that is invariant under projective automorphisms. Assume $Z(K_0(\mathbb{P}^n)) \subset \mathbb{Q} + i\mathbb{Q} $. Let $\mathcal{M}$ be a coherent sheaf on $\mathbb{P}^n \times \mathbb{A}^1$, flat over $\mathbb{A}^1$, whose fibers $M_t$ at $t\neq 0$ are projective transforms of a fixed nonzero sheaf $M$. If the special fiber $M_0$ has HN phases in $[v,u]$, where $u - v < 1$, then so does $M$.
\end{lemma}
\begin{proof}
    Choose $a< v \leq u < b$ with $b - a < 1$, such that the rays of phases $a$ and $b$ have rational slope. Consider the two hearts
    \[
    \mathscr{A}_a = \mathscr{P}((a,a+1]),
    \qquad
    \mathscr{A}_b = \mathscr{P}((b - 1,b]).
    \]
    Their intersection is $\mathscr{P}((a,b])$, and $M_0$ belongs to both. Both hearts are noetherian by \cite[Proposition 5.0.1]{AP06SheavesOftStructures}.

    Apply the open-heart property \cite[Proposition 3.3.2]{AP06SheavesOftStructures} to $M_0\in\mathscr{A}_a$, which gives a neighborhood $U_a$ of zero on which $\mathcal{M}$ lies in the relative heart. Curve generic flatness \cite[Lemma 3.5.2]{AP06SheavesOftStructures} then gives a dense open subset $W_a\subset U_a$ whose derived fibers belong to $\mathscr{A}_a$. Repeating the argument for $\mathscr A_b$ gives $W_b\subset U_b$. Choose $t\in W_a\cap W_b$ with $t\ne0$. Ordinary flatness identifies the derived fiber with $M_t = \mathcal{M}|_t$, so
    \[
    M_t \in \mathscr{A}_a \cap \mathscr{A}_b = \mathscr{P}((a,b]).
    \]
    Projective invariance implies $\phi^\pm_\tau(M)=\phi^\pm_\tau(M_t)$, so $M\in\mathscr P((a,b])$. Letting $a$ increase to $v$ and $b$ decrease to $u$ yields $\phi^-(M) \geq v$ and $\phi^+(M)\leq u$.
\end{proof}

The finite projection reduction used below appeared in Simpson's proof of boundedness of semistable sheaves \cite[Theorem $1.1$]{Simpson1994ModuliI}.

\begin{proposition}[phase bounds for pure sheaves on projective space]\label{proposition: phase bounds for pure sheaves on projective space}
    For a bounded family of nonzero pure $d$-dimensional sheaves $\mathscr{B}$ on $\mathbb{P}^n$, with $1\leq d\leq n$,
    \begin{equation}
        \phi^{\pm}_{\tau_{r,\lambda}}(M) = \psi_d(\lambda) + O(r^{-1})
        \qquad
        (M\in\mathscr{B}, \lambda\in V).
    \end{equation}
    Zero-dimensional sheaves have phase $1$ exactly.
\end{proposition}

\begin{proof}
    The case of $d = n$ is settled in \Cref{lemma:Support on linear subspace}. Assume $1\leq d < n$, and let $Z$ be the scheme-theoretic support of $M$. Choose a linear projection $f: Z\rightarrow \mathbb{P}^d$ with center disjoint from $Z$. A general center of dimension $n - d -1$ has this property. The projection is finite: the closure of a projection fiber is a linear space whose hyperplane at infinity lies in the center. Any positive dimensional fiber would meet that hyperplane, contradicting avoidance of the center. Thus the proper map is quasi-finite, hence finite.

    Put $V = f_* M$. It is torsion-free on $\mathbb{P}^d$. Indeed, its torsion subsheaf is preserved by the finite algebra $f_* \mathscr{O}_Z$ and hence corresponds to a coherent subsheaf of $M$, whose support has dimension less than $d$, contradicting purity unless it is zero. 
    
    The sheaves $V$ form a bounded family as $M$ and the projection vary. Indeed, choose one $m$ for which all $M$ are $m$-regular. Since $f$ is finite and $f^*\mathcal O(1)=\mathcal O_Z(1)$, we have $H^j(\mathbb P^d,V(k))=H^j(\mathbb P^n,M(k))$ for all $j,k$. Thus the $V$ are $m$-regular and have finitely many Hilbert polynomials. Their evaluation maps give a common quotient presentation, so they are bounded by \cite[Theorem 2.1]{Grothendieck95SGAIV}.
    
    Choose coordinates $[x:y]$ such that $f$ is $[x:y] \mapsto [x]$ and the center is $x = 0$. The map
    \begin{equation}
        h: Z\times \mathbb{A}^1 \rightarrow \mathbb{P}^n\times \mathbb{A}^1,
        \qquad
        ([x:y],t) \mapsto ([x:ty], t)
    \end{equation}
    is everywhere defined because $x\neq 0$ on $Z$. It is proper over $\mathbb{A}^1$, and each fiber map is finite: for $t\neq 0$ it is a projective transform of the inclusion $Z\hookrightarrow \mathbb{P}^n$; while for $t = 0$ it is $i\circ f$, where $i:\mathbb{P}^d \hookrightarrow \mathbb{P}^n$ is $y = 0$. Therefore $h$ is finite.

    Set
    \[
    \mathcal{M} = h_* \mathrm{pr}^*_Z M.
    \]
    This sheaf is flat over $\mathbb{A}^1$. And we could describe its fibers as above:
    \begin{equation}
        M_t \cong (g_t)_* M (t\neq 0),
        \qquad
        M_0 \cong i_*V,
        \qquad
        g_t([x:y]) =[x: ty].
    \end{equation}
    The special fiber's phase is bounded uniformly by \Cref{lemma:Support on linear subspace}. To apply \Cref{Lemma:Phase estimate from special fiber to general fiber}, we first suppose $r^j\lambda_j\in\mathbb Q+i\mathbb Q$ for every $j$, then the descended charge is rational on the numerical lattice and we have the desired bound for $M$. For arbitrary $\lambda$, fix $r$ and approximate $\lambda$ by such coefficient vectors in the slightly larger deformation neighborhood. The constants for the special-fiber estimate are uniform there. Continuity of the descended HN phases and of $\psi_d$ therefore passes the same bound to $\lambda$. This explicitly removes the rationality restriction.
\end{proof}

\begin{proof}[Proof of \Cref{Theorem: Dimension determines limit phase}]
    For $F\in\mathscr B$, put $G=\pi_*F$. Finite pushforward preserves purity and boundedness. The counit for $\pi^*\dashv\pi_*$ and the unit for $\pi_*\dashv\pi^!$ yield coherent maps
    \[
    \pi^*G\twoheadrightarrow F\hookrightarrow\pi^!G
    =\pi^*(G(n+1)).
    \]
    Then the estimate follows from \Cref{proposition: phase bounds for pure sheaves on projective space} and \Cref{lemma: Phase Control By Morphism}. 
    
    For the mass bound, by boundedness we have $|Z_{r,\lambda}(F)|=O_{\mathscr B}(r^d)$. Apply \Cref{lemma:mass-central-charge} with centre $\psi_d(\lambda)$ and $\varepsilon=M/r$. For large $r$,
    \[
    m_{r,\lambda}(F)\leq 2|Z_{r,\lambda}(F)|=O_{\mathscr B}(r^d).
    \]
\end{proof}

\subsection{Coefficient estimates}

We use the coefficient deformations from \Cref{Lemma:Deformations with uniform support constant} in two ways. First, we vary $\lambda_d$. Small mass and integrality force the relevant coefficient to vanish for the semistable object and every HN factor that can occur along the path. Its phase remains fixed, while the limiting phase of a pure $d$-dimensional sheaf moves; this gives the Hom vanishings below. For the coefficient estimate, $\lambda_d$ instead stays fixed and the other coordinates vary, allowing us to test the non-leading coefficients of a subobject in the moving heart.

\begin{lemma}\label{Lemma:Small-mass remains semistable in deformation}
Fix $1\leq d\leq n$ and a compact path in $V$ starting at $\zeta$ and changing only $\lambda_d$. Suppose its slicing distance from $\mathscr P_r$ is less than a fixed $a<1/2$, uniformly for $r\geq1$. There is $\varepsilon_d>0$, independent of $r$, such that every $\sigma_r$-semistable object $T$ with $m_r(T)<\varepsilon_dr^d$ remains semistable along the path, with its original phase.
\end{lemma}
\begin{proof}
    We first show that the charge of $T$ remains unchanged, and then that the charge of every object occurring as an HN factor anywhere on the path is independent of the deformation parameter. These are separate steps: constancy of $Z(T)$ alone would not prevent $T$ from becoming unstable.

    By \eqref{equation: Support property on coefficients} and integrality of $b_j$, sufficiently small $ \varepsilon_d$ forces $b_j(T) = 0$ for all $j \geq d$. In particular $Z(T)$ is unchanged. Its HN factors at any point along the path have phases in $[\phi_r(T) - a, \phi_r(T) + a]$. By \Cref{lemma:mass-central-charge}, we have
    \begin{equation}
        m_{r,\lambda}(T) \leq \sec(\pi a ) m_r (T).
    \end{equation}
    For every HN factor $U$ occurring anywhere on the path,
    \[
    r^{j} |b_j(U)| \leq C |Z_{r,\lambda}(U)| \leq C \sec(\pi a ) m_r(T) \leq C \sec(\pi a )\varepsilon_d r^d.
    \]
    Choose $\varepsilon_d$ so that $C\sec(\pi a )\varepsilon_d < 1$. When $j \geq d$ and $r\geq 1$, the inequality above gives
    \[
    |b_j(U)| \leq C \sec(\pi a )\varepsilon_d r^{d - j} < 1.
    \]
    Integrality forces $b_j(U) = 0$. Thus the charge of every actual HN factor also remains constant all along the path. For fixed $r, T, \zeta$, there are therefore finitely many coefficient tuples and constant nonzero charges. Each charge has at most one phase lift in the displayed interval of length $2a<1$. The continuous functions $\phi^+_{r,\lambda}(T)$ and $\phi^-_{r,\lambda}(T)$ take values in a finite set, and hence must be constant. They agree at the starting point, which proves the assertion.
\end{proof}

\begin{corollary}\label{Corollary:No Hom with small mass}
    Let $\mathscr{B}$ be a bounded family of pure $d$-dimensional sheaves, for $1\leq d \leq n$. There are $\rho,\varepsilon > 0$ and $r_0$ such that, for $F\in\mathscr{B}$, $r\geq r_0$, and any $\sigma_r$-semistable object $T$ with $m_r(T) < \varepsilon r^d$,
    \begin{equation*}
        \phi_r(T) < \phi_d + \rho \Longrightarrow \Hom (F, T) = 0,
    \end{equation*}
    \begin{equation*}
        \phi_r(T) > \phi_d - \rho \Longrightarrow \Hom (T,F) = 0.
    \end{equation*}
\end{corollary}
\begin{proof}
    Inside the neighborhood of \Cref{Lemma:Deformations with uniform support constant}, choose two small paths changing only $\lambda_d$ and ending at phases $\psi_d = \phi_d + 2\rho$ and $\psi_d = \phi_d - 2\rho$, respectively. Pick $\varepsilon$ to work for both paths in \Cref{Lemma:Small-mass remains semistable in deformation}. On the upward path the object $T$ remains semistable at its original phase, while \Cref{Theorem: Dimension determines limit phase}
    gives constants $M, r_0$ such that
    \[
    |\phi^{\pm}_{r,\lambda}(F) - \psi_d(\lambda) | \leq M/r \qquad (F\in\mathscr{B}, r\geq r_0).
    \]
    We might pick $r_0$ big enough such that 
    \[
    |\phi^{\pm}_{r,\lambda}(F) - \psi_d(\lambda) | < \rho \Longrightarrow \phi^-(F) \geq \phi_d +2\rho - \rho = \phi_d +\rho.
    \]
    Phase comparison then proves the first vanishing. On the downward path, the argument is the same.
\end{proof}

\begin{lemma}[Coefficient estimate]\label{Lemma: Coefficient estimate}
    Fix a bounded family $\mathscr B$ of nonzero pure $d$-dimensional sheaves, with $1\leq d\leq n$, a constant $c>0$, and a constant $M$ as in \Cref{Theorem: Dimension determines limit phase} for the family $\mathscr{B}$. There are constants $C,R>0$ such that every $F\in\mathscr B$ lies in $\mathscr A_r$ for $r\geq R$, and every nonzero subobject $J\hookrightarrow F$ in $\mathscr A_r$ satisfying
    \[
    \phi^-_r(J)\geq c,
    \qquad 
    \varphi_r(J) = \pi^{-1} \arg Z_r(J)\geq\phi_d-M/r
    \]
    satisfies
    \[
    m_r(J)\leq Cr^d,\qquad b_j(J)=0\ (j>d),
    \qquad
    0<b_d(J)\leq C,
    \qquad
    \sum_{j<d}r^j|b_j(J)|\leq Cr^{d-1}.
    \]
    All constants are independent of $F$, $J$, and $r$.
\end{lemma}
\begin{proof}
The case $J=F$ follows from boundedness and \Cref{Theorem: Dimension determines limit phase}. Assume henceforth that $Q=F/J\neq0$, where the quotient is taken in $\mathscr A_r$. Set
\[
\eta\coloneqq\tfrac18\min\{c,\phi_d,1-\phi_d\},
\qquad
\kappa\coloneqq\sin(\pi\eta)>0.
\]
Choose $0<\varepsilon\leq1$ so that the closed coefficient polydisc
\[
V'=\left\{\lambda_0=-1,\quad\lambda_d=\zeta_d,\quad
 |\lambda_j-\zeta_j|\leq\varepsilon
 \quad(1\leq j\leq n,\ j\neq d)\right\}
\]
is contained in $V$ and its slicings have distance at most $\eta$ from $\mathscr P_r$, uniformly in $r$. Increase $R$ so that every $F\in\mathscr B$ lies in $\mathscr A_r$ and $M/r\leq\eta$ for $r\geq R$.

The original exact sequence and \Cref{lemma: Phase Control By Morphism} give
\[
\phi^-_r(J)\geq c,\quad
\phi^+_r(J)\leq\phi_d+\eta,\quad
\phi^-_r(Q)\geq\phi_d-\eta,\quad
\phi^+_r(Q)\leq1.
\]
The corresponding bounds for $F$ are supplied by \Cref{Theorem: Dimension determines limit phase}. Each endpoint moves by at most $\eta$ under deformation. Since $c\geq8\eta$ and $8\eta\leq\phi_d\leq1-8\eta$, all three terms lie in
$\mathscr P_{r,\lambda}((2\eta,1+2\eta])$. Hence their original distinguished triangle is a short exact sequence
\begin{equation}\label{Equation: Short exact sequence in Lemma: Coefficient estimate for subobjects}
0\longrightarrow J\longrightarrow F\longrightarrow Q\longrightarrow0
\end{equation}
in this rotated heart, for every $\lambda\in V'$.

Put $\ell(z)=\im(e^{-2\pi i\eta}z)$. This functional is nonnegative on the charge of every object of the rotated heart. If $U$ is a deformed HN factor of $J$, then
\[
\eta\leq c-3\eta\leq \phi_{r,\lambda}(U)-2\eta\leq1-\eta.
\]
Here the upper bound follows already from $\phi^+_r(J)\leq1$ and the slicing-distance bound. Thus
\[
\ell(Z_{r,\lambda}(U))
 =|Z_{r,\lambda}(U)|\sin\bigl(\pi(\phi_{r,\lambda}(U)-2\eta)\bigr)
 \geq\kappa|Z_{r,\lambda}(U)|.
\]
Summing over the HN factors and using \eqref{Equation: Short exact sequence in Lemma: Coefficient estimate for subobjects}, we obtain
\begin{equation}\label{Equation: mass bound for subobject in Coefficient estimate for subobjects}
\kappa m_{r,\lambda}(J)
 \leq\ell(Z_{r,\lambda}(J))
 \leq\ell(Z_{r,\lambda}(F))
 \leq|Z_{r,\lambda}(F)|
 \leq K r^d,
\end{equation}
where $K$ depends only on the bounded family and the fixed neighborhood. Write $B=K/\kappa$; in particular, $m_{r,\lambda}(J)\leq Br^d$ uniformly on $V'$.

Define
\[
H_\lambda(P)=\im\bigl(e^{-\pi i\phi_d}Z_{r,\lambda}(P)\bigr).
\]
Because $J$ remains a subobject of $F$ in the rotated heart and $\lambda_d=\zeta_d$, we have
\[
\phi^+_{r,\lambda}(J)
 \leq\phi^+_{r,\lambda}(F)
 \leq\phi_d+M/r.
\]
Every HN phase $t$ of $J$ satisfies $t>0$. If $t-\phi_d<0$, then $t-\phi_d\in(-1,0)$ and its contribution to $H_\lambda(J)$ is nonpositive. Otherwise its sine is at most $\pi M/r$. Consequently, with $A=\pi MB$,
\[
H_\lambda(J)\leq A r^{d-1}\qquad(\lambda\in V').
\]
At $\zeta$, the total-charge phase lies between
$\phi_d-M/r$ and $\phi_d+M/r$ by hypothesis and the HN upper bound. Using $|Z_r(J)|\leq Br^d$ and increasing $R$ if necessary gives
\[
-A r^{d-1}\leq H_\zeta(J)\leq A r^{d-1}.
\]

We now vary only the free coordinates of $V'$:
\[
\lambda_0=-1,\qquad\lambda_d=\zeta_d,\qquad
\lambda_j=\zeta_j+i\varepsilon\zeta_d\operatorname{sgn}(b_j(J))
 \quad(1\leq j\leq n,\ j\neq d).
\]
The straight segment joining $\zeta$ to this point stays in $V'$. Since $|\zeta_d|=1$,
\[
\varepsilon\sum_{\substack{1\leq j\leq n\\j\neq d}}
 r^j|b_j(J)|
 =H_\lambda(J)-H_\zeta(J)
 \leq2A r^{d-1}.
\]
Set $S=2A/\varepsilon$. To estimate the fixed constant coefficient, use $\zeta_0=-1$ and compute
\[
H_\zeta(J)=\sin(\pi\phi_d)b_0(J)
 +\sum_{\substack{1\leq j\leq n\\j\neq d}}
 \sin\bigl(\pi(\phi_j-\phi_d)\bigr)r^j b_j(J).
\]
Since $0<\phi_d<1$, this implies
\[
|b_0(J)|\leq\frac{A+S}{\sin(\pi\phi_d)}r^{d-1}.
\]
Combining the two bounds yields a constant $C_1$, independent of $F,J,r$, such that
\begin{equation}\label{equation:single estimate for sums of coefficients in Lemma: Coefficient estimate for subobjects}
\sum_{j\neq d}r^j|b_j(J)|\leq C_1r^{d-1}.
\end{equation}
For $j>d$ and $r\geq1$, this gives $|b_j(J)|\leq C_1/r$. Increasing $R$ beyond $C_1$, integrality forces $b_j(J)=0$ for all $j>d$.

It remains to prove positivity and boundedness of $b_d(J)$. Write
\[
Z_r(J)=\zeta_d b_d(J)r^d+E_r,
\qquad |E_r|\leq C_1r^{d-1}.
\]
If the integer $b_d(J)$ were negative, then
\[
\im Z_r(J)\leq-\sin(\pi\phi_d)r^d+C_1r^{d-1}<0
\]
for sufficiently large $r$, contrary to $J\in\mathscr A_r$. Thus $b_d(J)\geq0$. The normalization $|\zeta_d|=1$ and the triangle inequality give
\[
|b_d(J)|r^d\leq|Z_r(J)|+|E_r|
 \leq Br^d+C_1r^{d-1},
\]
so $b_d(J)$ has a uniform upper bound.

Suppose finally that $b_d(J)=0$. Then $|Z_r(J)|\leq C_1r^{d-1}$. For sufficiently large $r$ the HN phases of $J$ lie in
\[
[c,\phi_d+M/r]
 \subset\bigl(0,(1+\phi_d)/2\bigr].
\]
Projection to the middle ray $(1+\phi_d)/4$ and \Cref{lemma:mass-central-charge} therefore give
\[
m_r(J)\leq
 \frac{C_1}{\cos(\pi(1+\phi_d)/4)}r^{d-1}.
\]
Apply \Cref{Corollary: Small mass forces support drop} with $k=d-1$ and this fixed mass constant. After one further uniform increase of $R$, every coherent cohomology sheaf of $J$ has support dimension at most $d-1$. By \Cref{Lemma:Geometric Amplitude}, $J\in\mathrm D^{\leq0}_{\Coh}(C^n)$, and standard truncation gives
\[
\Hom(J,F)=\Hom(\mathscr H^0(J),F)=0,
\]
the last equality following from purity of $F$. This contradicts the nonzero monomorphism $J\hookrightarrow F$. Thus $b_d(J)>0$.

Taking $C$ larger than the mass, coefficient-sum and leading-coefficient constants proves every assertion. All thresholds and constants depend only on $\mathscr B,c,M$ and the fixed geometric and deformation data, and not on $F,J,r$.
\end{proof}

\section{Gieseker implies Limit}

All constants in this section are uniform over the indicated bounded family and the subobjects in question. Restrictions to divisors are derived and regarded as objects on $C^n$ by pushforward.

The main step is to show that an exact sequence concentrated near $\phi_d$ is eventually a saturated coherent sequence. Applied to HN and Jordan--H\"older filtrations, this will reduce stability to the phase comparison of \Cref{lemma:asymptotic-sheaf-phase-comparison}.

\begin{theorem}\label{theorem: Narrow inclusions become saturated and coherent}
Let $\mathscr B$ be a bounded family of pure $d$-dimensional sheaves, where $1\leq d\leq n$, and let $M>0$ be the constant given by \Cref{Theorem: Dimension determines limit phase}. There is $R>0$ such that, for $E\in\mathscr B$ and $r\geq R$, every exact sequence
    \[
    0\longrightarrow A\longrightarrow E\longrightarrow B\longrightarrow0
    \]
in $\mathscr A_r$ whose nonzero terms have HN phases in $[\phi_d-M/r,\phi_d+M/r]$ is a coherent exact sequence with $A\subset E$ saturated. As $E$, the sequence, and $r\geq R$ vary, the sheaves $A$ and $B$ form bounded families.
\end{theorem}

Here saturation means that the quotient is zero or pure of dimension $d$. We prove the theorem by bounding the coherent image of $A\to E$, then showing that $A$ equals its saturation. 

\begin{proposition}\label{proposition: coherent image estimate}
    Fix a bounded family $\mathscr{B}$ of pure $d$-dimensional sheaves, $d\geq 1$. Fix $0 < c \leq \phi_n/2$ and $M>0$ the constant given by \Cref{Theorem: Dimension determines limit phase}. There are constants $K,R_0$ such that the following holds for every $r\geq R_0$. Suppose
    \begin{equation}\label{exact sequence in  coherent image estimate}
        0\rightarrow A \rightarrow F \rightarrow B \rightarrow 0,
        \qquad
        F\in\mathscr{B},
        \qquad
        A\neq 0,
    \end{equation}
    is exact in $\mathscr{A}_r$ and 
    \begin{equation}
        \phi^-_r(A) \geq c, 
        \qquad
        \varphi_r(A) \geq \phi_d - M/r,
        \qquad
        \left|\phi^{\pm}_r(B) - \phi_d \right| \leq M/r,
    \end{equation}
    when $B \neq 0$. For $I = \text{im}(\mathscr{H}^0(A) \rightarrow F)$, one then has
    \begin{equation}
        b_d(I) = b_d(A) > 0,
        \qquad
        b_{d - 1}(I)  \geq -K.
    \end{equation}
\end{proposition}

\begin{proof}
    We induct on $d$. All constants in the proof are fixed by $\mathscr{B}$, $c$, $M$ and $\zeta$. If $B = 0$, the assertion follows trivially from boundedness of $F$. There are some preparations before we start: by \Cref{Lemma: Coefficient estimate}, there is one constant $C$ such that
    \begin{equation}\label{equation: coefficients control preparation}
        m_r(A) \leq C r^d,
        \qquad
        b_j(A) = 0 \quad (j > d),
        \qquad
        0< b_d(A) \leq C,
        \qquad
        |b_{j}(A)| \leq C r^{d - j - 1} \quad (j < d).
    \end{equation}
    Next, fix a sufficiently large degree $a$ such that $aH$ is very ample. We choose $D \in |aH|$ smooth and regular on all standard cohomology sheaves of $A$, $F$, $B$, also regular on $F/I$, and so that $F|_D$ is pure. For purity, we also require $D$ to be regular on $\mathcal E xt_{C^n}^q(F,\omega_{C^n})$ for all $q$; \cite[Lemma 1.1.12 and Corollary 1.1.14 (ii)]{HL10GeometryOfModuliOfSheaves} then apply. Only finitely many additional open conditions are imposed. The restricted $\{F|_D \}$ form a bounded pure family, as they are obtained from a finite-type family and a fixed linear system. 
    
    By \Cref{lemma: fixed twists}, we might pick $r$ large enough so that $\eqref{exact sequence in  coherent image estimate}$ remain exact in $\mathscr{A}_r$ under twist by $\mathcal{O}_X(-D)$. We could enlarge $r$ so that the bounded coherent sequence
    \[
    0\rightarrow F(-D) \rightarrow F \rightarrow F|_D \rightarrow 0
    \]
    is exact in $\mathscr{A}_r$.     
    \[
    \begin{tikzcd}[
    ]
    & {}
    & 0 \arrow[d]
    & {}
    & {}
     \\
    0 \arrow[r]
    & A(-D) \arrow[r] \arrow[d, hook]
    & F(-D) \arrow[r] \arrow[d]
    & B(-D) \arrow[r] \arrow[d]
    & 0
    \\
    0 \arrow[r]
    & A \arrow[r] \arrow[d]
    & F \arrow[r] \arrow[d, two heads]
    & B \arrow[r]
    & 0
    \\
    & \left.A\right|_D \arrow[r]
    & \left.F\right|_D
    & {}
    & {}
    \end{tikzcd}
    \]
    The snake lemma then gives a short exact sequence in $\mathscr{A}_r$:
    \begin{equation}\label{snake seqeuence in hyperplane induction}
        0\rightarrow R \rightarrow A|_D \rightarrow J \rightarrow 0,
    \end{equation}
    where $J\subset F|_D$, $R = \ker (B(-D) \rightarrow B)$. And we have the phase bounds by \Cref{lemma: Phase Control By Morphism}:
    \[
    \phi^+_r(R) \leq \phi_d + M/r,
    \qquad
    \phi^-_r(J) \geq c,
    \qquad
    \phi^+_r(J) \leq \phi_{d-  1} + M/r.
    \]
    The first follows from $R\subset B(-D)$, the second from the successive quotients $A \twoheadrightarrow A|_D \twoheadrightarrow J$, and the third from $J \subset F|_D$.
    
    \paragraph{\textbf{Case $d = 1$}} Numerical restriction bounds $b_0(A)$ and $b_1(A)$. It suffices to prove in this case both $A$ and $B$ are coherent sheaves.  Since $Z_r(A|_D) = -ab_1(A) < 0$, $R$ and $J$ can have only phase-one factors. But $\phi^+_r(R) \leq \phi^+_r(B(-D)) < 1$, so $R = 0$. We have the full exact sequence in $\mathscr{A}_r$:
    \[
    0\rightarrow A|_D \rightarrow F|_D\rightarrow B|_D \rightarrow 0.
    \]
    \Cref{Corollary: Small mass forces support drop} puts all standard cohomology of $A,B$ in support dimension at most one, hence all standard cohomology of $A|_D$ and $B|_D$ is zero-dimensional. Their Postnikov towers coincide with their HN filtrations, and $A|_D$, $B|_D \in \mathscr{A}_r$ force their negative cohomology to vanish. Regularity now implies the negative cohomology of $A$, $B$ is zero-dimensional; the same argument then forces them to vanish. Hence both $A$ and $B$ are coherent sheaves, and $I = A$.

    \paragraph{\textbf{Case $d > 1$}} Now suppose $d \geq 2$. The main estimate below bounds $m_r(R)$ by $O(r^{d-2})$, so every standard cohomology sheaf of $R$ has support dimension at most $d-2$. A controlled enlargement by \Cref{lemma: controlled enlargement} then supplies the quotient phase bound required by induction. Its coherent image differs from $I|_D$ by a sheaf of dimension at most $d-2$ and bounded leading coefficient, which suffices to recover the lower bound on $b_{d-1}(I)$.
    
    Since $J$ has HN phases in $[c, (1+\phi_{d - 1})/2]\subset(0,1)$, let 
    \[
    \kappa = \min \{ \sin(\pi c), \sin (\pi(1+\phi_{d - 1})/2)\},
    \] 
    when $J\neq 0$, we have
    \[
    \kappa m_r(J) \leq \im Z_r(J) \leq \im (Z_r(F|_D)) \leq |Z_r(F|_D)| \leq C_1 r^{d - 1},
    \]
    where $C_1$ is a constant depending only on the bounded family $\{F|_D\}$. Take $L_{d - 1} = C_1/\kappa$, we conclude that
    \[
    m_r(J) \leq L_{d - 1}r^{d - 1}.
    \]
    Riemann--Roch and \eqref{equation: coefficients control preparation} give a uniform constant $D_{d - 1}$ with
    \begin{equation}\label{equation: coefficient estimate for A|_D}
        \left|   Z_r(A|_D)    -adb_d(A)\zeta_{d - 1}r^{d - 1}   \right| \leq D_{d - 1} r^{d - 2}.
    \end{equation}

    The key observation here is: $R$ can only have support dimension at most $d - 2$ when $r\gg 0$. Indeed, \eqref{equation: coefficient estimate for A|_D} yields:
    \begin{equation}\label{equation: Im estimate for A|_D}
        \left| \im( e^{-\pi i \phi_{d - 1}} Z_r(A|_D) ) \right| \leq D_{d - 1}r^{d - 2}.
    \end{equation}
    From the phase bound $\phi^+_r(J) \leq \phi_{d-  1} + M/r$ and $m_r(J) \leq L_{d - 1} r^{d - 1}$,  by \Cref{lemma:charge-phase-mass relation}
    \begin{equation}\label{equation: Im estimate for J}
         \im( e^{-\pi i \phi_{d - 1}} Z_r(J) )  \leq m_r(J) \cdot \left( \frac{\pi M}{ r}\right) \leq \pi M L_{d - 1}r^{d - 2}.
    \end{equation}
    We might pick $r$ large enough such that $\phi^+_r(R)\leq \phi_d + M/r \leq \frac{\phi_{d - 1}+ \phi_d}{2}$, then by \Cref{lemma:charge-phase-mass relation}
    \begin{equation}\label{equation: Im estimate for R}
        \im(e^{-\pi i  \phi_{d- 1}} Z_r(R))\leq m_r(R)\cdot \sin\left(-\pi \frac{\phi_{d - 1} - \phi_d}{2}\right) = -\sin\left(\frac{\pi - \theta}{2}\right)m_r(R).
    \end{equation}
    Combining \eqref{equation: Im estimate for A|_D}, \eqref{equation: Im estimate for J}, \eqref{equation: Im estimate for R}, we obtain
    \begin{equation}
        m_r(R) \leq L_R r^{d - 2},
        \qquad
        L_R \coloneqq \frac{\pi ML_{d - 1} + D_{d - 1}}{\sin((\pi - \theta)/2)} 
    \end{equation}
    We might enlarge $r$ so that by \Cref{Corollary: Small mass forces support drop}, $b_{d -1 }(R) = 0$. Hence
    \begin{equation}
        b_{d - 1} (J) = b_{d - 1}(A|_D) = adb_d(A) > 0, 
        \qquad
        \left| \varphi_{r}(J) - \phi_{d - 1}\right| \leq M/r.
    \end{equation}
    The estimate for the total-charge phase of $J$ does not yet bound the maximal HN phase of its quotient in $F|_D$. By \Cref{lemma: controlled enlargement}, there is $J\subset A' \subset F|_D$ with $B'\coloneqq F|_D/ A'$, such that 
    \[
    0 \rightarrow A' \rightarrow F|_D \rightarrow B' \rightarrow 0
    \] 
    satisfies the $(d - 1)$-induction hypotheses, and
    \[
    b_{d - 1}(A') = b_{d - 1}(J) = adb_d(A).
    \]
    Write $I_J = \text{im}(\mathscr{H}^0(J) \rightarrow F|_D)$, $I' = \text{im}(\mathscr{H}^0(A') \rightarrow F|_D)$, then $I_J = I|_D$ by regularity, and from the same lemma, we have $I|_D = I_J \subset I'$ with
    \begin{equation}
        \dim \mathrm{Supp}(I'/I|_D) \leq d - 2,
        \qquad
        0\leq b_{d - 2}(I'/I|_D)  \leq \Delta,
    \end{equation}
    for a constant $\Delta$. From the induction step, there is a constant $K_{d - 1}$ such that
    \begin{equation}
        b_{d - 1}(I') = b_{d - 1}(A'), 
        \qquad
        b_{d - 2}(I') \geq -K_{d - 1}.
    \end{equation}
    We have
    \begin{align*}
        b_{d - 1}(I|_D) &= ad b_d(I) > 0, \\
        b_{d - 2}(I|_D) &= a(d - 1)b_{d - 1}(I) - a^2\binom{d}{2}b_d(I).
    \end{align*}
    Setting $K\coloneqq \frac{1}{a(d - 1)}\left(K_{d - 1} + \Delta \right)$, we obtain:
    \begin{equation*}
        b_d(I) = b_d(A) > 0,
        \qquad
        b_{d - 1}(I)  \geq -K. \qedhere
    \end{equation*}
\end{proof}

The following controlled enlargement supplies the quotient phase bound needed for induction on the support dimension.

\begin{lemma}\label{lemma: controlled enlargement}
    Fix a bounded family $\mathscr{B}$ of pure $(d-1)$-dimensional sheaves, $d\geq 2$. Fix $0 < c \leq \phi_n/2$ and $M>0$ the constant given by \Cref{Theorem: Dimension determines limit phase}. There are constants $M'\geq M$, $L, \Delta ,R$ depending only on these data and the fixed angle, with the following property. 
    
    For every $r\geq R$, suppose $E\in \mathscr{B}$, $0\neq J \subset E$  is a subobject in $\mathscr{A}_r$, satisfying
    \begin{equation}
        \phi^-_r(J) \geq c, 
        \qquad
        \varphi_r(J) \geq \phi_{d - 1} - M/r,
        \qquad
        \left| \phi^\pm_r(E) - \phi_{d - 1}  \right| \leq M/r.
    \end{equation}
    Then there is an enlargement
    \begin{equation}
        J\subset A'\subset E,
        \qquad
        B' \coloneqq E/ A',
    \end{equation}
    such that
    \begin{equation}
        \phi^-_r(A') \geq c,
        \qquad
        \varphi_r(A') \geq \phi_{d - 1} - M/r,
        \qquad
        |\phi^\pm_r(B') - \phi_{d - 1} |\leq M'/r
    \end{equation}
    when $B' \neq 0$, and
    \begin{equation}
        b_{d - 1}(A'/J) = 0, 
        \qquad
        m_r(A'/J) \leq L r^{d - 2}.
    \end{equation}
    Moreover, the coherent images
    \[
    I_J = \text{im}(\mathscr{H}^0(J) \rightarrow E),
    \qquad
    I' = \text{im}(\mathscr{H}^0(A') \rightarrow E),
    \]
    satisfy 
    \begin{equation}
        I_J \subset I',
        \qquad
        \dim \mathrm{Supp}(I'/I_J) \leq d - 2,
        \qquad
        0\leq b_{d - 2}(I'/I_J)  \leq \Delta.
    \end{equation}
\end{lemma}
\begin{proof}
    We inspect the HN factors of $V=E/J$. A factor with positive leading coefficient already has phase within $O(r^{-1})$ of $\phi_{d-1}$. The factors above a suitable threshold therefore have zero leading coefficient. Absorbing this HN prefix into $J$ controls the remaining quotient. If $V = 0$, take $A' = E$ and all assertions follow. Otherwise let 
    \[
    0= V_0 \subset V_1 \subset \cdots \subset V_\ell = V
    \]
    be its HN filtration. For each $j$, let $P_j \subset E$ be the inverse image of $V_j$, thus $P_0 = J$. Also,
    \[
    \phi^-_r(V_j) \geq \phi^-_r(V) \geq \phi^-_r(E) \geq \phi_{d - 1} - M/r.
    \]
    The exact sequences 
    \[
    0\rightarrow J \rightarrow P_j \rightarrow V_j \rightarrow 0
    \]
    then yield
    \begin{align*}
        \phi^-_r(P_j) &\geq \min \{ \phi^-_r(J), \phi^-_r(V_j) \} \geq \min \{ \phi^-_r(J), \phi^-_r(E)\} \geq c, \\
        \varphi_r(P_j)  &\geq \min \{\varphi_r(J),  \phi^-_r(V_j)  \}\geq \phi_{d - 1}  -M/r.
    \end{align*}
    Thus \Cref{Lemma: Coefficient estimate} applies uniformly for every $P_i$, with one common constant $C_0$. In particular,
    \begin{equation*}
        Z_r(P_j) = \zeta_{d- 1} r^{d - 1} b_{d - 1}(P_j) + e(P_j),
        \qquad
        |e(P_j)| \leq C_0 r^{ d - 2}.
    \end{equation*}
    Hence for the HN factors $V_j/V_{j - 1}$, we have
    \begin{equation*}
        Z_r(V_{j}/V_{j - 1}) = \zeta_{d-  1}r^{d - 1}b_{d - 1}(V_j/V_{j - 1}) + e_j,
        \qquad
        |e_j| \leq 2C_0 r^{d - 2}.
    \end{equation*}
    For $r > 2C_0/\sin(\pi\phi_{d - 1})$, integrality of $b_{d - 1}$ and non-negativity of $\im Z_r(V_j/V_{j - 1})$ imply \[b_{d -1}(V_j/V_{j - 1})\geq 0.\] If $b_{d -1}(V_j/V_{j - 1})> 0$, for $r \geq 4C_0$, the elementary inequality $ |\arg(1+w)| \leq 2 |w|$ for $|w| \leq 1/2$ implies
    \begin{equation}\label{inequality: positive leading term implies phase narrowness}
        |\phi_r(V_j/V_{j - 1})-\phi_{d - 1}|\leq C_1/r,\qquad C_1=4C_0/\pi.
    \end{equation}
    Choose $M'=\max\{M,C_1\}$. Now set
    \[
    h \coloneqq \max\{j \mid \phi_r(V_j/V_{j - 1}) > \phi_{d - 1} + M'/r \},
    \]
    with $h = 0$ if there are no such factors, and take $A' \coloneqq P_{h}$, then $A'/J = V_h$. We claim $A'$ satisfies all the assertions. From \eqref{inequality: positive leading term implies phase narrowness}, every factor $V_j/V_{j - 1}$ with $j\leq h$ has $b_{d - 1}=0$, so $b_{d - 1}(V_h)=0$. Every remaining factor has phase at most $\phi_{d - 1}+M'/r$ and at least $\phi_{d - 1}-M/r$. This proves the asserted phase bounds for $B'=V/V_h$. The bounds for $A'$ are the already established bounds for $P_h$.
    
    Next we bound the mass of $V_h$. The charge estimates for $P_h$ and $P_0 = J$ yield
    \[
    |Z_r(V_h)| \leq 2 C_0 r^{d - 2}.
    \]
    All HN phases of $V_h$ lie in $[\phi_{d - 1}, 1]$, take $L\coloneqq \frac{2C_0}{\cos (\pi(1 - \phi_{d - 1})/2)}$. Using \Cref{lemma:mass-central-charge} with the middle ray $(\phi_{d - 1} + 1)/2$ yields
    \[
    m_r(V_h) \leq \frac{1}{\cos (\pi(1 - \phi_{d - 1})/2)} |Z_r(V_h)| \leq  L r^{d - 2}.
    \]
    Finally, since objects in $\mathscr{A}_r$ have no positive coherent cohomology, the exact sequence in $\mathscr{A}_r$
    \begin{equation}
        0\rightarrow J \rightarrow A' \rightarrow A'/J \rightarrow 0
    \end{equation}
    gives
    \begin{equation}
        \mathscr{H}^0(J) \rightarrow \mathscr{H}^0(A') \rightarrow \mathscr{H}^0(A'/J) \rightarrow 0.
    \end{equation}
    Consequently $I_J \subset I'$, and $I'/I_J$ is a quotient of $\mathscr{H}^0(A'/J)$. The mass estimate of $A'/J$ and \Cref{Corollary: Small mass forces support drop} now finish all the remaining assertions.
\end{proof}

We next use Grothendieck's boundedness lemma in the form that allows the ambient sheaf to vary.

\begin{lemma}[{\cite[Proposition 1.2(i), Lemma $2.6$]{Grothendieck95SGAIV}}]\label{lemma: Grothendieck's boundedness lemma}
    Fix a bounded family of pure $d$-dimensional sheaves $\mathscr{B}$, $d\geq 1$, and a constant $c$. For any $E\in\mathscr{B}$, the saturated subsheaves $0\neq S \subset E$ satisfying $b_{d - 1}(S) \geq -c$ form a bounded family. Their pure quotients $E/S$ are bounded as well.
\end{lemma}

Combined with the coherent image estimate above, now we could prove \Cref{theorem: Narrow inclusions become saturated and coherent}.

\begin{proof}[Proof of \Cref{theorem: Narrow inclusions become saturated and coherent}]
    If $A=0$ or $B=0$, the conclusion is immediate; assume both are nonzero. Let $I$ be the coherent image of $A\rightarrow E$, $S\coloneqq \mathrm{Sat}_E(I)$ and $Q \coloneqq E/S$. By the coherent image estimate \Cref{proposition: coherent image estimate},
    \[
    b_d(S) = b_d(I) = b_d(A) > 0,
    \qquad
    b_{d - 1}(S) \geq b_{d - 1}(I) \geq -C.
    \]
    Grothendieck boundedness \Cref{lemma: Grothendieck's boundedness lemma} then bounds $S$, $Q$, and their coherent sequence with $E$ is uniformly exact in $\mathscr{A}_r$ for large $r$ by \Cref{CoherentSheavesInLargeVolumeHeart}.

    The map $A \rightarrow E$ factors through $S$, and $A \rightarrow S$ is also monic in $\mathscr{A}_r$. Its cone $C$ fits into
    \[
    0\rightarrow A \rightarrow S \rightarrow C \rightarrow 0,
    \qquad
    0\rightarrow C \rightarrow B \rightarrow Q \rightarrow 0.
    \]
    If $C \neq 0$, by the quotient and subobject phase bounds, we have
    \[
    \phi^-_r(S) \leq \phi^-_r (C) \leq \phi^+_r (C) \leq \phi^+_r(B).
    \]
    Thus $|\phi^{\pm}_r(C) - \phi_d | \leq M/r$. Since $b_d(S) = b_d(A)$, by \Cref{Lemma: Coefficient estimate}, we have $|Z_r(C)| = O(r^{d- 1})$, and the phase concentration implies $m_r(C) = O(r^{d - 1})$. The last HN factor $T_r$ of $C$ is a nonzero semistable quotient of $S$, with $m_r(T_r) = O(r^{d - 1})$ and $\phi_r(T_r) = \phi_d + O(r^{-1})$, which leads to a contradiction by \Cref{Corollary:No Hom with small mass} for the bounded pure family $S$. Therefore $C = 0$ and $A = S$ as subobjects of $E$, proving every assertion.
\end{proof}

\begin{corollary}\label{corollary:HN coefficients estimate}
For a bounded family $\mathscr B$ of nonzero pure $d$-dimensional sheaves, $d\geq1$, all HN factors for sufficiently large $r$ are pure $d$-dimensional sheaves and form a bounded family. The same holds for Jordan--H\"older factors of semistable members. In particular, all their coefficients are uniformly bounded, and each factor $T_i$ of $E$ satisfies $0<b_d(T_i)\leq b_d(E)$.
\end{corollary}
\begin{proof}
The phase estimate and \Cref{theorem: Narrow inclusions become saturated and coherent} apply to every filtration prefix $E_i\subset E$. Thus the quotients $Q_i=E/E_i$ form a bounded family of pure $d$-dimensional sheaves or zero. Each factor $T_i=E_i/E_{i-1}$ is a saturated subsheaf of $Q_{i-1}$, with quotient $Q_i$. Moreover,
\[
 b_{d-1}(T_i)=b_{d-1}(Q_{i-1})-b_{d-1}(Q_i)
\]
is uniformly bounded below. Grothendieck boundedness bounds the factors. Purity and additivity give the asserted bounds on $b_d$.
\end{proof}

\begin{theorem}[Uniform comparison on $C^n$]\label{theorem: uniform Gieseker implies limit}
 Let $\mathscr B$ be a bounded family of pure $d$-dimensional sheaves on $C^n$, where $1\leq d\leq n$. There is $R>0$ such that, for every $E\in\mathscr B$ and every $r\geq R$, Gieseker semistability of $E$ implies $\sigma_r$-semistability, and Gieseker stability implies $\sigma_r$-stability.
\end{theorem}
\begin{proof}
 By \Cref{Theorem: Dimension determines limit phase}, all HN phases of $E\in\mathscr B$ lie in $[\phi_d-M/r,\phi_d+M/r]$ for uniform $M$ and large $r$. If $E$ is not $\sigma_r$-semistable, let $A\subset E$ be its first HN term. By \Cref{theorem: Narrow inclusions become saturated and coherent}, this is a proper saturated coherent subsheaf, and all such $A$ form a bounded family. After a further uniform increase of $R$, \Cref{lemma:asymptotic-sheaf-phase-comparison} gives
 \[
 p_H(A)\leq p_H(E)\quad\Longrightarrow\quad
 \varphi_r(A)\leq\varphi_r(E).
 \]
 This contradicts the defining inequality for the first HN term, and proves the semistability assertion.

 If $E$ is Gieseker stable but strictly $\sigma_r$-semistable, choose the first nonzero proper term $A\subset E$ of a Jordan--H\"older filtration. Both it and its quotient have the phase of $E$, so the same theorem and uniform comparison apply. Gieseker stability gives $\varphi_r(A)<\varphi_r(E)$, contradicting equality of these phases. Taking the maximum of the two uniform thresholds proves the claim.
\end{proof}

\begin{remark}
 In dimension zero, every nonzero sheaf is Gieseker semistable and $\sigma_r$-semistable of phase one. In either theory it is stable precisely when it is a skyscraper sheaf. Thus the comparison also holds for bounded families of zero-dimensional sheaves.
\end{remark}

\subsection{Passage to a smooth projective variety}\label{subsec:general-comparison}

Return to the fixed embedding $i_X:X\hookrightarrow\mathbb P^N$, with $H=i_X^*\mathcal O_{\mathbb P^N}(1)$, and fix $\theta\in(\pi-\pi/N,\pi)$. We work at radii satisfying $r\sin\theta>\alpha_{i_X}$, so the induced stability conditions exist by \Cref{StabilityConditionOnSmoothProjectiveVarieties}. For the finite flat map $\pi:C^N\to\mathbb P^N$ of \Cref{StabilityConditionOnPn}, write
\[
 \Phi=\pi^*i_{X*}:\Db(X)\longrightarrow\Db(C^N).
\]
The induced slicing and the normalization of the charge give, for every nonzero $E$,
\begin{equation}\label{eq:comparison-transfer}
 \phi^\pm_{r,\theta,X}(E)=\phi^\pm_{2r,\theta,C^N}(\Phi E),
 \qquad
 Z_{2r,\theta,C^N}(\Phi E)=(\deg\pi)Z_{r,\theta,X}(E).
\end{equation}
Here the product radius is $2r$, since $\pi^*\mathcal O_{\mathbb P^N}(1)=\mathcal O_{C^N}(2H)$, and $\deg\pi=2^NN!$.

The functor $\Phi$ is exact and conservative for standard cohomology and exact between the induced hearts. It sends an HN filtration to an HN filtration, proving the phase identity in \eqref{eq:comparison-transfer}, and preserves nonzero proper subobjects. It also preserves bounded families of pure sheaves and their support dimension. For purity, finite pushforward preserves support dimension, and $\pi_*\pi^*G=G\otimes\pi_*\mathcal O_{C^N}$ is locally a finite direct sum of $G$; hence a lower-dimensional subsheaf of $\pi^*G$ would contradict purity of $G$.

\begin{theorem}[Uniform comparison on $X$]\label{Thm:UniformComparisonOnX}
Let $\mathscr B$ be a bounded family of nonzero pure coherent sheaves on $X$. There is $R>0$ such that, for every $E\in\mathscr B$ and $r\geq R$, the sheaf $E$ lies in $\mathscr A_{r,\theta,X}$, and $H$-Gieseker (semi)stability of $E$ implies $\sigma_{r,\theta,X}$-(semi)stability.
\end{theorem}
\begin{proof}
We may fix the support dimension $d$ and then take the maximum of the finitely many thresholds. If $d=0$, a finite-length filtration by skyscraper sheaves is exact in the induced heart. Geometricity makes every factor stable of phase $1$, so every such sheaf is semistable, and it is stable precisely when it has length one, in either theory. Assume $d\geq1$.

Apply \Cref{Theorem: Dimension determines limit phase} to the bounded pure family $\Phi\mathscr B$. By \eqref{eq:comparison-transfer}, there is a uniform $M$ such that all HN phases of every $E\in\mathscr B$ lie in $[\phi_d-M/r,\phi_d+M/r]$. This interval is contained in $(0,1)$ for large $r$, so $E$ lies in the induced heart.

If $E$ is unstable, let $A$ be its first HN subobject. If $E$ is semistable but not stable, let $A$ be a nonzero proper subobject of the same phase. The exact sequence
\[
 0\longrightarrow A\longrightarrow E\longrightarrow B\longrightarrow0
\]
and its image under $\Phi$ have all HN phases in the same interval. Apply \Cref{theorem: Narrow inclusions become saturated and coherent} on $C^N$ at radius $2r$, with constant $2M$. It follows that $\Phi A\subset\Phi E$ is a saturated coherent subsheaf and that all such $\Phi A$ form a bounded family. Since $\Phi$ detects standard cohomology, $A\subset E$ is a nonzero proper coherent subsheaf on $X$.

To compare phases uniformly, only finitely many charge polynomials are needed. Comparing coefficients in \eqref{eq:comparison-transfer} gives
\begin{equation}\label{eq:transfer-coefficients}
 b_j(\Phi A)=(\deg\pi)2^{-j}\tilde a_j(A).
\end{equation}
Boundedness of $\Phi A$ therefore gives finitely many Todd-twisted Hilbert polynomials for $A$, while boundedness of $\mathscr B$ gives the same for $E$. The polynomial argument of \Cref{lemma:asymptotic-sheaf-phase-comparison} and the triangular relation in \Cref{ToddTwistDoesNotAffectGieseker} give one threshold after which the comparisons $<,=,>$ between $\varphi_{r,\theta,X}(A)$ and $\varphi_{r,\theta,X}(E)$ agree with those between $p_H(A)$ and $p_H(E)$ for all these pairs.

A first HN subobject has strictly larger phase, so it would contradict Gieseker semistability of $E$. Thus a Gieseker semistable sheaf is Bridgeland semistable. If $E$ is Gieseker stable, a proper same-phase subobject would instead give $p_H(A)=p_H(E)$, again a contradiction. This proves the stable assertion as well.
\end{proof}

\begin{proof}[Proof of \Cref{Thm:MainComparison}]
A fixed numerical class determines the Hilbert polynomial, so the pure $H$-Gieseker semistable sheaves of class $v$ form a bounded family by \cite[Theorem $1.1$]{Simpson1994ModuliI}; see also \cite[Theorem 4.2]{Langer2004SemistablePositiveChar}, \cite[Theorem 3.3.7]{HL10GeometryOfModuliOfSheaves}. Apply \Cref{Thm:UniformComparisonOnX} to obtain one threshold for both semistability and stability. The implication `Limit implies Gieseker' follows from \Cref{LimitsImpliesGiesekerGeneral}, applied to each object along the entire ray tail.
\end{proof}

\begin{remark}[Scope of the argument]\label{Rem:ScopeOfComparison}
The product geometry enters in the uniform deformation and mass estimates of Section~\ref{sec:uniform-estimates}, and in the phase estimates of Section~\ref{section: Pure sheaf phase estimates and coefficients control}. Once these inputs and eventual membership of bounded sheaf families in the hearts are available, the coefficient, coherent-image, and boundedness arguments use general properties of coherent sheaves and stability conditions. On $X$ these inputs need not be reproved: the functor $\Phi$ transfers precisely the phase and coherence statements required above. In particular, the proof does not require $\Phi$ to preserve Gieseker stability or stable objects.
\end{remark}

\bibliography{BridgelandGiesekerILimitCase}
\bibliographystyle{alpha}

\end{document}